\documentclass[12pt]{article}

\usepackage{graphicx}
\usepackage[francais]{babel}
\usepackage{amsfonts}
\usepackage{pstricks}
\usepackage{amsmath}
\usepackage{amssymb}
\usepackage{amsthm}
\usepackage{empheq}
\usepackage{indentfirst}
\input{amssym.def}
\input{amssym}
\allowdisplaybreaks

\renewcommand{\theequation}{\arabic{section}.\arabic{equation}}

\newcommand{\R}{\mathbb{R}}

\newcommand{\pa}{\mathbb{\partial}}
\newcommand{\eps}{\varepsilon}

\def\dis{\displaystyle}

\newtheorem{lemma}{\bf Lemma}[section]
\newtheorem{theorem}{\bf Theorem}[section]
\newtheorem{prop}{\bf Proposition}[section]
\newtheorem{corol}{Corollary}[section]
\newtheorem{defin}{Definition}[section]
\newtheorem{remark}{\bf Remark}[section]

\def\md{\mathrm{d}}

\title{Long-time behaviors and stability of entropy solutions for
linearly degenerate hyperbolic systems of rich type}

\author{Yue-Jun PENG$^{\;1,2}$ and Yong-Fu YANG$^{\;3,4}$}

\date{}

\begin{document}
\maketitle

\vspace{-3mm}

\begin{center}
{\small $^1$Clermont Universit\'e, Universit\'e Blaise Pascal, 63000 Clermont-Ferrand, France \\[2mm]
$^2$CNRS-UMR 6620, Laboratoire de Math\'ematiques, 63171 Aubi\`ere cedex, France \\[2mm]
$^3$LCP, Institute of Applied Physics and Computational Mathematics, Beijing, 100088, China \\[2mm]
$^4$Department of Mathematics, College of Sciences, Hohai University, Nanjing 210098, China \\[2mm]
Email : peng@math.univ-bpclermont.fr, \hspace{2mm}yongfuyang@sina.com}
\end{center}

% --------------------------------------------------------------------------

\vspace{5mm}

\begin{center}
\begin{minipage}{14cm}
        {\bf Abstract.} {\small We show that in one space dimension, a linearly
degenerate hyperbolic system of rich type admits exact traveling wave solutions
if the initial data are Riemann type outside of a space interval. In a particular
case of the system including physical models, we prove the convergence of entropy
solutions to traveling waves in the $L^1$ norm as the time goes to infinity. The
traveling waves are determined explicitly in terms of the initial data and the
system. We also obtain the stability of entropy solutions in $L^1$.}
\end{minipage}
\end{center}

\vspace{5mm}


%$1%%%%%%%%%%%%%%%%%%%%%%%%%%%%%%%%%%%%%%%%%%%%%%%%%%%%%%%%%%%%%%%%%%%%%%

\section{Introduction}
\newcounter{intro}
\renewcommand{\theequation}{\thesection.\theintro}

For first order quasilinear hyperbolic systems of conservation laws, it is well
known that, generically speaking, classical solutions exist only locally in time
and singularities may appear in a finite time. Therefore, we are often led to seek
global entropy solutions (see \cite{Jo90, La57, Ser96} and the references
therein). In some special cases, however, the global existence of classical solutions
can be obtained. In particular, these situations occur for linearly degenerate systems
provided that the initial data are sufficiently small and decay at infinity
(see \cite{CH04, Li94, LZK94, Lin04}).

As to the Cauchy problem in one space dimension, the existence of global entropy
solutions was proved by Glimm \cite{Gl65} when the initial data are small in
the total variation norm. The long-time behavior of the entropy solution was
investigated by Glimm and Lax \cite{GL70} for systems of two conservation
laws having genuinely nonlinear characteristics. In \cite{Liu77} Liu studied
this problem for a general hyperbolic system of conservation laws when the initial
data agree with Riemann data outside of a large space interval. He showed
that waves of genuinely nonlinear characteristic fields tend to the solution of
the corresponding Riemann
problem and waves of linearly degenerate characteristic fields tend to traveling
waves, both at algebraic rates. The convergence is in the BV norm.
When the initial data have compact support, the solution also converges in the
$L^1$ norm to the linear superposition of traveling waves and $N$-waves.

If the system is strictly hyperbolic and weakly linearly degenerate (see \cite{Li94}
for definition), Kong and Yang \cite{KY03} also discussed the long-time behavior
of global classical solutions to its Cauchy problem with small and decaying initial
data, provided that the source term satisfies a so-called matching condition.
Recently, Liu and Zhou \cite{LZ07} considered the similar long-time behavior of
classical solutions for diagonal systems when the initial data are large. They
applied their results to the generalized extremal surface equations. All these
results indicate that the classical solutions of a linearly degenerate hyperbolic
system approach some $C^1$ traveling wave solutions. Unfortunately, all these results
are qualitative on the limits, from which we do not know what these traveling waves
are.

Stability is another important problem for hyperbolic systems of conservation
laws. In particular, it implies the uniqueness of solutions.  A complete theory
of the uniqueness
and $L^1$ stability of entropy solutions to the scalar conservation law was
established by Kruzkov \cite{Kr70}. See also \cite{BP98} and the references therein
by a kinetic method. For hyperbolic systems of conservation laws, we refer to the
work of Bressan \cite{Br00} for the existence, uniqueness and stability of entropy
solutions in $L^\infty(\R) \cap BV(\R)$ with small initial data.

In a previous work, the authors of this paper considered the generalized extremal
surface equations (see \cite{Br02, KSZ06}) and showed that its Lipschitz solutions
are equivalent to entropy
solutions in $L^{\infty}(\R)$ of a non-strictly hyperbolic system of conservation laws
\cite{PY10}. They obtained an explicit representation formula and the uniqueness
of the entropy solutions to the Cauchy problem. Based on these results, they further
obtained the convergence and convergence rates of the entropy solutions to traveling
waves in the $L^1(\R)$ norm. The traveling waves are determined explicitly in terms
of the initial data and the system. Moreover, when initial data agree
with Riemann data outside of a finite space interval, the entropy solutions become
explicit traveling waves after a finite time. Finally, they also got $L^1$ stability
of the entropy solutions.

In this paper, we are going to extend results of \cite{PY10} to linearly degenerate
hyperbolic systems of rich type. The notion of rich system was introduced by Serre
\cite{Ser89} and Tsar\"ev \cite{Ts85, Ts90} independently. It is also called
semi-hamiltonian system by Tsar\"ev. For this kind of linearly degenerate system,
the global existence of entropy
solutions in $L^\infty$ was established by Chen \cite{Ch92} in the strictly hyperbolic
case and in \cite{LPR09} in the non-strictly hyperbolic case. The uniqueness of entropy
solutions has not been solved yet. However, it is solved in a particular case where
$N_i=N$ (see the definition below) for all $i$. This case contains many interesting
physical models, such as the generalized extremal surface equations, the Born-Infeld
system and augmented Born-Infeld system. All these models are non-strictly hyperbolic
systems.

Our study on the long-time behavior and stability of entropy solutions in $L^1$
is carried out in this special case. The main results are stated in Theorems
\ref{T3.1} and \ref{T4.1}, which are available for initial data in $L^1 \cap L^\infty$
and also in $C^k$ with integer $k \geq 1$. The traveling wave solution after a finite
time holds in a
more general case, provided the initial data are Riemann type outside of a finite
space interval (see Theorem \ref{T2.1}). Below are preliminary descriptions on the
rich system and related results.

Consider the Cauchy problem for quasilinear hyperbolic systems of diagonal form~:
 \stepcounter{intro}
\begin{equation}
\label{diag}
    \pa_t \,w_i+\lambda_i(w)\, \pa_x \,w_i=0,
\quad \forall \, i \in K_n, \quad t > 0, \quad x \in \R
\end{equation}
with initial data
\stepcounter{intro}
\begin{equation}
\label{diag-init}
t=0~:\quad w=w^0(x),\quad x \in \R.
\end{equation}
Here $K_n = \{1,\cdots,n\}$, $w=(w_1,\cdots,w_n)^{\top}$ and $w^0=(w^0_1,\cdots,w^0_n)^{\top}$, in which $\top$ means transpose. The eigenvalues $\lambda_i(w)\; (i \in K_n)$ of the system are real valued smooth
functions defined on an open domain of $\R^n$. We suppose that each of them has
a constant multiplicity, namely, on the domain under consideration we have
\stepcounter{intro}
\begin{align}
\nonumber
    \mu_1(w)\overset{\text{def}}{=}\lambda_1(w)=\cdots=\lambda_{r_1}(w)
& < \mu_2(w) \overset{\text{def}}{=}\lambda_{r_1+1}(w)
=\cdots=\lambda_{r_2}(w)\\ \label{hyper}
    &< \cdots < \mu_s(w) \overset{\text{def}}{=}\lambda_{r_{s-1}+1}(w)
=\cdots=\lambda_{r_s}(w),
\end{align}
where $r_1,\;r_2,\ldots,\;r_s$ are constants and
$$1 \leq r_1 < r_2< \cdots < r_s=n.$$
For system (\ref{diag}) and $p \in K_s$, the characteristic
$\mu_p(w)$ with constant multiplicity $r_p-r_{p-1} \; (r_0=0)$ is
linearly degenerate if and only if $\mu_p(w)$ is independent
of $w_{r_{p-1}+1},\cdots,w_{r_p}$~:
$$  \frac{\pa\mu_p(w)}{\pa w_l}\equiv 0, \quad (l=r_{p-1}+1,\cdots,r_p). $$
If all characteristics $\lambda_i(w) \, (i \in K_n)$ are linearly degenerate,
system (\ref{diag}) is said to be linearly degenerate. Let
$\Lambda(w)=\text{diag}\{\lambda_1(w),\cdots,\lambda_n(w)\}$. A pair of functions
$(E(w),F(w))$ is an entropy-entropy flux pair of (\ref{diag}) if
$F'(w)=E'(w)\Lambda(w)$, i.e.,
$$\frac{\pa F(w)}{\pa w_j}=\lambda_j(w)\frac{\pa E(w)}{\pa w_j},
\quad \forall \, j \in K_n.$$
It is well known that $(E(w),F(w))$ is an entropy-entropy flux pair if and only
if it satisfies the additional conservation law (see \cite{La57,Ser96})
$$      \pa_t E(w)+\pa_x F(w)=0.    $$

Recall that the $i$th characteristic is rich if, on the domain under consideration,
for all $j\,,k\in K_n$ such that $\lambda_j(w)\neq \lambda_i(w)$ and
$\lambda_k(w)\neq \lambda_i(w)$, we have
$$   \dis \frac{\pa}{\pa w_j}\Big(\dfrac{\frac{\pa
\lambda_i(w)}{\pa w_k}}{\lambda_k(w)-\lambda_i(w)}\Big)
=\frac{\pa}{\pa w_k}\Big(\dfrac{\frac{\pa
\lambda_i(w)}{\pa w_j}}{\lambda_j(w)-\lambda_i(w)}\Big).   $$
This definition is equivalent to the existence of a smooth function
$N_i(w)>0$ such that
\stepcounter{intro}
\begin{equation}
\label{richness1}
\big(\lambda_j(w)-\lambda_i(w)\big)\frac{\pa N_i(w)}{\pa w_j}
=N_i(w)\frac{\lambda_i(w)}{\pa w_j},\quad \forall\, j \in K_n,
\quad \lambda_j(w)\neq \lambda_i(w).
\end{equation}
When $\lambda_i(w)$ is linearly degenerate, (\ref{richness1}) is equivalent to
$$ \big(\lambda_j(w)-\lambda_i(w)\big)\frac{\pa N_i(w)}{\pa w_j}
=N_i(w)\frac{\lambda_i(w)}{\pa w_j},\quad \forall\, j \in K_n.  $$
System (\ref{diag}) is said to be rich if all characteristics
$\lambda_i(w)\; (i \in K_n)$ are rich. We refer to \cite{Ser96,Sev94,Ts90}
for the definition and properties of the rich system.

For any given smooth solution $w=w(t,x)$ to the linearly degenerate strictly
hyperbolic rich system (\ref{diag}), from \cite{Ser96} we have
\stepcounter{intro}
\begin{equation}
\label{richness3}
\pa_t \,N_i(w)+ \pa_x \,\big(N_i(w)\lambda_i(w)\big)=0, \quad i \in K_n.
\end{equation}
Therefore, $(N_i(w),N_i(w)\lambda_i(w))$ is an entropy-entropy flux
pair of system (\ref{diag}). More generally, we have~:
\stepcounter{intro}
\begin{equation}
\label{entropy}
\pa_t\,\big(N_i(w)g_i(w_i)\big)+ \pa_x\,\big(N_i(w)\lambda_i(w)g_i(w_i)\big)
=0\quad \text{for\, any\, smooth\, function}\, g_i.
\end{equation}
Thus, $N_i(w)g_i(w_i)\; (i \in K_n)$ stand for $n$ independent families of
entropies and each entropy of the linearly degenerate rich system (\ref{diag})
is a linear combination of $N_i(w)g_i(w_i)\; (i \in K_n)$ (see \cite{Ser89,Ser96}).
In particular, taking $g_i(s)=s$ in (\ref{entropy}) gives~:
\stepcounter{intro}
\begin{equation}
\label{entropy1}
    \pa_t \,\big(N_i(w)w_i\big)+ \pa_x \,\big(N_i(w)\lambda_i(w)w_i\big)=0.
\end{equation}

\begin{defin}
\label{Def1}
A function $w\in L^{\infty}(\R^+\times \R)$ is an entropy solution to a linearly
degenerate hyperbolic rich system of diagonal form (\ref{diag}) if the entropy
equalities (\ref{richness3}) and (\ref{entropy1}) are satisfied in the sense of
distribution for all $i\in K_n$.
\end{defin}

Let $w\in L^{\infty}(\R^+\times \R)$ be an entropy solution. For each $i\in K_n$,
as a compatibility condition, (\ref{richness3}) implies that there is a function
$y_i=Y_i(t,x)$ such that
\stepcounter{intro}
\begin{equation}
\label{EL}
    \md y_i=N_i(w)\,\md x - (N_i\lambda_i)(w)\,\md t.
\end{equation}
Since
$$ \frac{\pa Y_i}{\pa x}|_{t=0}=N_i(w(t,x))|_{t=0}=N_i(w^0(x)),  $$
it is natural to define~:
\stepcounter{intro}
\begin{equation}
\label{EL-initial}
Y_i(0,x)=Y_i^0(x)\stackrel{\text{def}}{=}\int_0^x\,N_i(w^0(\xi))\,\md \xi.
\end{equation}
For any given $w\in L^{\infty}(\R^+\times \R)$, (\ref{EL}) and (\ref{EL-initial})
determine a unique Lipschitz and strictly increasing function $y_i=Y_i(t,x)$ for
all $t \geq 0$. As a consequence, from Proposition 2 and Theorem 3 in \cite{LPR09},
we have

\begin{prop}
\label{P1.1}
Suppose that $w^0 \in L^{\infty}(\R)$ and system (\ref{diag}) is linearly
degenerate and rich. Then the Cauchy problem (\ref{diag})-(\ref{diag-init})
admits a global entropy solution $w\in L^{\infty}(\R^+\times \R)$ in the sense of
Definition \ref{Def1}. Moreover,
\stepcounter{intro}
\begin{equation}
\label{exp-sol}
         w_i(t,x)=w_i^0(X_i^0(Y_i(t,x))) \quad (i \in K_n),
\end{equation}
where $X_i^0=(Y_i^0)^{-1}$ is the inverse function of $Y_i^0(\cdot)$. Moreover, if
$w^0 \in L^{\infty}(\R) \cap C^k(\R)$ for some integer $k \geq 1$, then
$w \in L^{\infty}(\R^+\times \R) \cap C^k(\R^+\times \R)$.
\end{prop}

\begin{remark}
\label{R1.1}
     Proposition \ref{P1.1} only establishes the existence of the entropy
solutions. Up to our knowledge, so far no uniqueness results
are available in this general case of the rich system.
\end{remark}

Next we consider the special case that $N_i=N$ for all $i \in K_n$. Then
(\ref{richness3}) becomes~:
$$  \pa_t \,N(w)+ \pa_x \,\big(N(w)\lambda_i(w)\big)=0, \quad i\in K_n. $$
From \cite{LPR09}, we know that for each $i \in K_n$, there are a constant
$\tilde{\lambda}_i$ and a function $M$ independent of $i$ such that
\stepcounter{intro}
\begin{equation}
\label{5.1}
     N(w)\lambda_i(w)-\tilde{\lambda}_i=M(w).
\end{equation}
Hence,
$$    \pa_t \,N(w)+ \pa_x \,M(w)=0. $$
This scalar equation allows to make a change of coordinates
$(t,x) \longmapsto (s,z)$ with (see \cite{Wa87, Pe07})
\stepcounter{intro}
\begin{equation}
\label{F5.1}
     s=t,\quad \md z =N\,\md x-M\,\md t.
\end{equation}
For $w^0\in L^{\infty}(\R)$, let $z=Z(t,x)$ be the unique Lipschitz function
satisfying (\ref{F5.1}) and
\stepcounter{intro}
\begin{equation}
\label{Z0}
Z(0,x)=Z^0(x)\overset{\text{def}}{=}\int_0^x\,N\big(w^0(\xi)\big)\,\md \xi.
\end{equation}
Since $x \longmapsto Z(t,x)$ is a strictly increasing function, we may define
its inverse by $X(t,\cdot) = Z^{-1}(t,\cdot)$ for all $t \geq 0$. The function
$x=X(t,x)$ satisfies
\stepcounter{intro}
\begin{equation}
\label{X1}
\md x =\frac{1}{N(\tilde{w}(t,z))}\,\md z
+\frac{M(\tilde{w}(t,z))}{N(\tilde{w}(t,z))}\,\md t,
\end{equation}
where
\stepcounter{intro}
\begin{equation}
\label{tdw}
    \tilde{w}(t,z)=\big(\tilde{w}_1^0(z-\tilde{\lambda}_1\,t),\cdots,
\tilde{w}_n^0(z-\tilde{\lambda}_n\,t)\big)^{\top}, \quad \tilde{w}_i^0 = w_i^0 \circ X^0.
\end{equation}
Together with $X(0,z)=X^0(z)$, (\ref{X1}) determines a unique Lipschitz function
for all $t\geq 0$. Then the unique entropy solution of the Cauchy problem
(\ref{diag})-(\ref{diag-init}) is given by
\stepcounter{intro}
\begin{equation}
\label{solu}
  w_i(t,x)=w_i^0\big(X^0(Z(t,x)-\tilde{\lambda}_i\,t)\big),
\quad \forall \, i \in K_n.
\end{equation}
Moreover, a straightforward computation using (\ref{5.1}) gives
\stepcounter{intro}
\begin{equation}
\label{X2}
     X^0(z) = \int_0^z \frac{1}{N(\tilde{w}^0(z))}\,\md z,
\quad X(t,z) = X^0(z) + \int_0^t \Big(\lambda_i(\tilde{w}(\tau,z))
- \frac{\tilde{\lambda}_i}{N(\tilde{w}(\tau,z))}\Big)\,\md \tau.
\end{equation}

Thus, Li \textit{et al.} \cite{LPR09} got the following result.

\begin{prop}
\label{P1.2}
Assume that system (\ref{diag}) is linearly degenerate and rich with
$N_i$ being independent of $i$ for all $i \in K_n$. Then the Cauchy
problem (\ref{diag})-(\ref{diag-init}) with the initial data $w^0\in
L^{\infty}(\R)$ admits a unique entropy solution $w \in
L^{\infty}(\R^+\times\R)$ given by expression (\ref{solu}), where
$X^0$ is the inverse function of $Z^0$ being defined by (\ref{Z0})
and $Z(t,\cdot)$ is the inverse function of the unique solution
$X(t,\cdot)$ to (\ref{X1}) and $X(0,z)=X^0(z)$.
\end{prop}

This paper is organized as follows. In the next section, we first
derive the explicit traveling waves after a finite time. Section 3
and Section 4 are devoted to the study of the long-time behavior and
$L^1$ stability of the entropy solutions in the case $N_i=N$ for all
$i$, respectively. Finally, we give two examples in Section 5.

%$2%%%%%%%%%%%%%%%%%%%%%%%%%%%%%%%%%%%%%%%%%%%%%%%%%%%%%%%%%%%%%%%%%%%%%%

\section{Explicit traveling waves after a finite time}
\newcounter{asymp}
\renewcommand{\theequation}{\thesection.\theasymp}

Throughout this section, we suppose that the initial data $w^0$ are constants
outside of a finite space interval $(-L,L)$ with $L> 0$. More precisely, let
\stepcounter{asymp}
\begin{equation}
\label{const}
w^0(x)=\bar{w}^\pm \stackrel{\text{def}}{=}
(\bar{w}_1^{\pm},\cdots,\bar{w}_n^{\pm})^{\top}, \quad \forall \pm x \geq L,
\end{equation}
where $\bar{w}_i^{\pm}$ are constants for all $i \in K_n$. In addition,
to simplify the presentation, we suppose temporarily that system (\ref{diag})
is strictly hyperbolic. Based on the explicit expression (\ref{exp-sol}), we
show that the entropy solution is exactly composed of traveling waves after
a finite time $t_L > 0$ to be determined later. In particular, each Riemann
invariant is a traveling wave after $t_L$.

For each $i\in K_n$, since $N_i(w) > 0$, (\ref{EL}) implies that
$x \longmapsto Y_i(t,x)$ is a strictly increasing function and then a bijection
for all $t \geq 0$. Let $X_i(t,\cdot)$ be the inverse function of $Y_i(t,\cdot)$
and
$$    X_i^{\pm}(t)=X_i(t,Y_i^0(\pm L)).     $$
It follows from (\ref{exp-sol}) and (\ref{const}) that
$$   w_i(t,x)=\bar{w}_i^+, \quad \mbox{if} \; \; X_i^0(Y_i(t,x)) \geq L,  $$
namely,
$$  w_i(t,x)=\bar{w}_i^+, \quad \forall\; x \geq X_i^+(t).  $$
Similarly,
$$  w_i(t,x)=\bar{w}_i^-, \quad \forall\; x \leq X_i^-(t).  $$
From the strict hyperbolicity, we have $X_j^\pm(t) < X_i^\pm(t)$ for all
$1 \leq j < i \leq n$. Then
\stepcounter{asymp}
\begin{equation}
\label{w_j^+}
   w_j(t,x)=\bar{w}_j^+, \quad \forall \; x \geq X_i^+(t), \quad \forall \; j \leq i
\end{equation}
and
\stepcounter{asymp}
\begin{equation}
\label{w_j^-}
   w_j(t,x)=\bar{w}_j^-, \quad \forall \; x \leq X_i^-(t), \quad \forall \; j \geq i.
\end{equation}
In particular, we obtain
\stepcounter{asymp}
\begin{equation}
\label{wXn^+}
    w(t,x) = \bar{w}^+, \quad \forall\; x \geq X_n^+(t).
\end{equation}
Similarly,
\stepcounter{asymp}
\begin{equation}
\label{wX1^-}
    w(t,x)=\bar{w}^-, \quad \forall\; x \leq X_1^-(t).
\end{equation}

\begin{remark}
\label{R2.1} $\; $

\noindent (1) For each $i\in K_n$, $x=X_i^\pm(t)$ is the $i$th
characteristic curve issued from the point $(0,\pm L)$, respectively.

Indeed, $X_i^\pm(0)= \pm L$ and it follows from (\ref{EL}) that
$$  \dis \frac{\md X_i^\pm(t)}{\md t} =\frac{\md X_i(t,Y_i^0(\pm L))}{\md t}
=\lambda_i(w(t,X_i^\pm)),   $$
where $w$ is given by (\ref{exp-sol}).

\vspace{3mm}

\noindent (2) The characteristics $x=X_1^-(t)$ and $x =X_n^+(t)$ are both
straight lines.

Indeed, (\ref{wXn^+}) and (\ref{wX1^-}) imply that
$$ \lambda_1(w(t,x)) = \lambda_1(\bar{w}^-) \; \; \text{if} \; \; x \leq X_1^-(t) $$
and
$$ \lambda_n(w(t,x)) = \lambda_n(\bar{w}^+) \; \; \text{if} \; \; x \geq X_n^+(t). $$
Obviously, $\lambda_1(\bar{w}^-)$ and $\lambda_n(\bar{w}^+)$ are constants, we then obtain
$$    X_1^-(t)=\lambda_1(\bar{w}^-)\,t-L,
\quad X_n^+(t)=\lambda_n(\bar{w}^+)\,t+L.   $$
\end{remark}

Since the system is linearly degenerate, for $1 \leq i < j \leq n$, two curves
$x=X_i^+(t)$ and $x=X_j^-(t)$ intersect only once. Let $t_{i, j} >0$ be their
intersection time. Then $X^+_i(t_{i, j}) = X^-_j(t_{i, j})$ for $1 \leq i < j \leq n$.
We denote
\stepcounter{asymp}
\begin{equation}
\label{t*}
    t_L = \max_{1 \leq i < j \leq n}\;t_{i, j}
\end{equation}
and
$$    D=\{(t,x)\,|\, x \in \R, \,t > t_L\}.   $$
Noting that the function $Y_i^0$ is strictly increasing, we deduce from the strict
hyperbolicity that
\stepcounter{asymp}
\begin{equation}
\label{orderX_i}
    X_1^-(t) < X_1^+(t) < X_2^-(t) < X_2^+(t) < \cdots < X_n^-(t) < X_n^+(t),
    \quad \forall \, t > t_L.
\end{equation}
Then $D$ can be divided into $2n+1$ sub-domains separated by the
characteristics $x= X_i^{\pm}(t)\;(i \in K_n)$. We denote these sub-domains
by $\dis D_0,\,D_i\;(i \in K_n),\,D_{n+j}\;(j \in K_{n-1})$ and $D_{2n}$,
respectively, where
\begin{align}
    \dis &D_0 = \{(t,x)\,|\, x \leq X_1^-(t),\, t > t_L\},\quad D_{2n}
= \{(t,x)\,|\, x > X_n^+(t),\, t > t_L\},\nonumber \\[2mm]
 \dis &D_i = \{(t,x)\,|\, X_i^-(t) < x \leq X_i^+(t),\, t > t_L\}
\quad (i \in K_n), \nonumber \\[2mm]
 \dis &D_{n+i} = \{(t,x)\,|\, X_i^+(t) < x \leq X_{i+1}^-(t),\, t
> t_L\} \quad (i \in K_{n-1}). \nonumber
\end{align}
From (\ref{wXn^+}) and (\ref{wX1^-}), we obtain
$$    w(t,x)=\bar{w}^-, \quad \forall \, (t,x) \in D_0   $$
and
$$    w(t,x)=\bar{w}^+, \quad \forall \, (t,x) \in D_{2n}.   $$
Moreover, (\ref{w_j^+}) and (\ref{w_j^-}) imply that
$$   w(t,x)=\bar{w}^i, \quad \forall \, (t,x) \in D_{n+i}
\; \;  (i \in K_{n-1}),  $$
where
$$  \bar{w}^i=(\bar{w}_1^+,\cdots,\bar{w}_i^+,\bar{w}_{i+1}^-,\cdots,\bar{w}_n^-)^{\top},
\quad \forall \, i \in K_n. $$

For $i \in K_n$, from (\ref{orderX_i}) we have
$$  X_{i-1}^+(t) < X_i^-(t) < x \leq X_{i+1}^-(t) < X_{i+1}^+(t),
\quad \forall \, (t,x) \in D_i, $$
which together with (\ref{w_j^+}) and (\ref{w_j^-}) gives
$$  w_j(t,x) = \bar{w}_j^+, \quad \forall \, (t,x) \in D_i,
\; \; \forall \, j < i  $$
and
$$  w_j(t,x) = \bar{w}_j^-, \quad \forall \, (t,x) \in D_i,
\; \; \forall \, j > i. $$
Since $\lambda_i$ is independent of $w_i$, we obtain
$$   \lambda_i(w(t,x)) = \lambda_i(\bar{w}^i), \quad \forall \, (t,x) \in D_i.  $$
Hence, equation (\ref{diag}) shows that $w_i$ is a traveling wave in $D_i$
with speed $\lambda_i(\bar{w}^i)$. In view of expression (\ref{exp-sol}),
we may write
\stepcounter{asymp}
\begin{equation}
\label{Di}
  w(t,x)=\big(\bar{w}_1^+,\cdots,\bar{w}_{i-1}^+,w_i^0(\phi_i(x-\lambda_i(\bar{w}^i) t))
,\bar{w}_{i+1}^-,\cdots,\bar{w}_n^-\big)^{\top}, \quad \mbox{if} \; \; (t,x) \in D_i,
\end{equation}
where $\phi_i\; (i \in K_n)$ are Lipschitz functions because $X^0_i$ and $Y_i$
are Lipschitzian. Remark that $X^-_i$
and $X^+_i$, $X^-_i$ and $X^+_j \, (i < j)$ do not intersect after time $t_L$.
The domains $D_i$ and curves $X^\pm_i$ are illustrated in Figure 1.

\vspace{5mm}
\begin{center}
\includegraphics{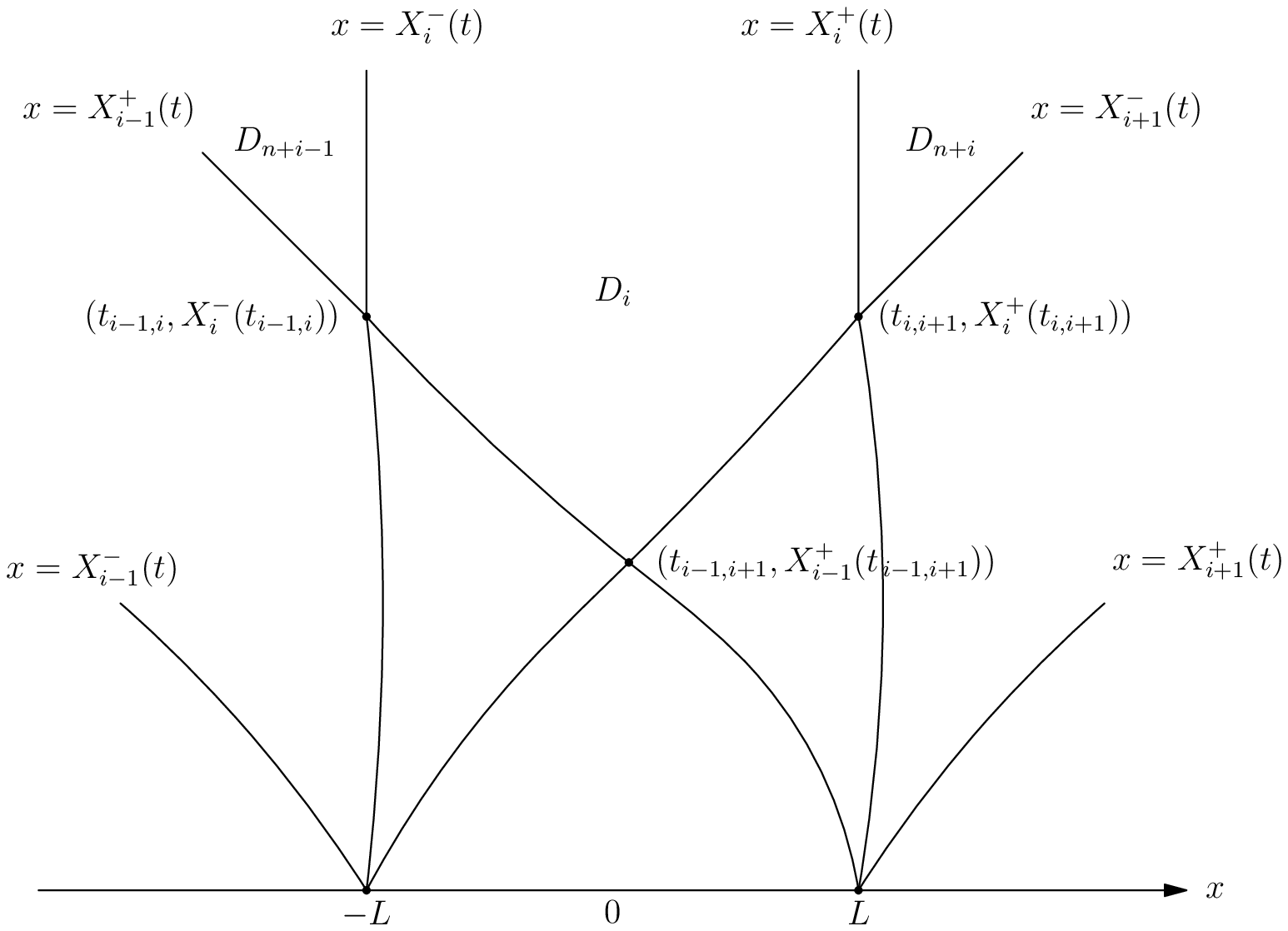}

\vspace{2mm} Figure 1. Traveling wave solutions after a finite time.
\end{center}

Thus, we get the main result of this section.

\begin{theorem}
\label{T2.1}
Let the assumptions of Proposition \ref{P1.1} hold. Assume furthermore that
(\ref{const}) holds. Then there exists a finite time $t_L>0$ defined by
(\ref{t*}) and Lipschitz functions $\phi_i \; (i \in K_n)$ such that
$$    w(t,x)= \left\{\begin{array}{llll}
(\bar{w}_1^-,\cdots,\bar{w}_n^-)^{\top}, \quad (t,x) \in D_0, \\[3mm]
(\bar{w}_1^+,\cdots,\bar{w}_{i-1}^+,w_i^0(\phi_i(\eta_i)),\bar{w}_{i+1}^-,
\cdots,\bar{w}_n^-)^{\top},\;(t,x) \in D_i, \; \; i \in K_n, \\[3mm]
    (\bar{w}_1^+,\cdots,\bar{w}_i^+,\bar{w}_{i+1}^-,\cdots,\bar{w}_n^-)^{\top},
\quad (t,x) \in D_{n+i}, \; \; i \in K_{n-1}, \\[2mm]
(\bar{w}_1^+,\cdots,\bar{w}_n^+)^{\top}, \quad (t,x) \in D_{2n},
\end{array}\right. $$
where $\eta_i=x-\lambda_i(\bar{w}^i) t$.
\end{theorem}

\begin{remark}
\label{R2.2}
Similar results are still valid under the assumption that system (\ref{diag})
has characteristics with constant multiplicity. Indeed, if (\ref{hyper}) holds,
then there are essentially $s$ distinct characteristics $\mu_p
\; (p \in K_s)$. Similarly to (\ref{t*}), we may define a finite time
still denoted by $t_L>0$. In this situation, $D$ is divided into
$2s+1$ sub-domains~:
$$  D_p = \{(t,x)\,|\, X_p^-(t) < x \leq X_p^+(t),\; t > t_L\},
\quad p \in K_s,    $$
$$  D_p=\{(t,x)\,|\, X_p^+(t) < x \leq X_{p+1}^-(t),\; t >
t_L\}, \quad p \in K_{s-1}. $$
On each domain $D_{s+p}$ we have
$$  w(t,x)=\bar{w}^{\,r_p}, \; \; (t,x) \in D_{s+p} $$
and
$$  w(t,x)=\Big(\bar{w}_1^+,\cdots,\bar{w}_{r_{p-1}}^+,w^0_{r_{p-1}+1}(\phi_p(\eta_p)),
\cdots,w^0_{r_p}(\phi_p(\eta_p)),\bar{w}_{r_p+1}^-,\cdots,\bar{w}_n^-\Big)^{\top}, \; \; (t,x) \in D_p, $$
where $\phi_p$ is a Lipschitz function and $\eta_p=x-\mu_p(\bar{w}^{\,r_p}) t$.
\end{remark}

%$3%%%%%%%%%%%%%%%%%%%%%%%%%%%%%%%%%%%%%%%%%%%%%%%%%%%%%%%%%%%%%%%%%%%%%%

\section{Long-time behaviors of entropy solutions}
\newcounter{long}
\renewcommand{\theequation}{\thesection.\thelong}

This section is devoted to the study of long-time behaviors of the global
entropy solutions in the case that $N_i=N$ for all $i \in K_n$. By means
of the explicit expression of the entropy solutions and the special
structure of the system, we prove that as $t$ goes to infinity, the entropy
solutions converge to traveling wave solutions, which are also of explicit
forms in terms of the initial data and the system. More precisely, we assume
that there is a constant vector $\bar{w}=(\bar{w}_1,\cdots,\bar{w}_n)^{\top}$
such that $w^0-\bar{w} \in L^1(\R) \cap L^{\infty}(\R)$ and
$$    \lim_{x\rightarrow \pm \infty} w^0(x)=\bar{w}.   $$
In particular, if $w^0-\bar{w}\in W^{1,1}(\R)$, the above limits hold.
We need also smallness conditions on the initial data in appropriate norms
according to the situation $\tilde{\lambda}_i=0$ or $\tilde{\lambda}_i \neq 0$
(see Lemmas 3.5-3.6). The convergence result holds in the sense of $L^1(\R)$
norm (see Theorem \ref{T3.1}).

Let us denote
$$    \bar{\lambda}_i=\lambda_i(\bar{w}_1,\cdots,
\bar{w}_{i-1},\bar{w}_{i+1},\cdots,\bar{w}_n), \quad \forall \,i \in K_n.  $$
The first two lemmas are consequences of the structure of the system.

\begin{lemma}
\label{L3.1}
    Under the assumptions of Proposition \ref{P1.2}, for each
$i \in K_n$ we have
\stepcounter{long}
\begin{eqnarray}
\label{eq0}
    X\big(t,Z^0(x)+\tilde{\lambda}_i\,t\big)-\bar{\lambda}_i\,t
    & \!\!=\!\! & x + \int_{Z^0(x)}^{Z^0(x)+\tilde{\lambda}_i\,t}
\Big(\frac{1}{N\big(\tilde{w}^0(\xi)\big)}-\frac{1}{N(\bar{w})}\Big)
\, \md\xi \nonumber \\[2mm]
    & \!\!\;\!\! & + \int_0^t \Big(\frac{M}{N}(\tilde{w}(\tau,Z^0(x)
+\tilde{\lambda}_i\,t))-\frac{M}{N}(\bar{w})\Big)\,\md \tau,
\end{eqnarray}
where $\tilde{w}$ and $\tilde{w}^0$ are defined by (\ref{tdw}).
\end{lemma}

\noindent {\bf Proof.} From (\ref{X1}) we have
\begin{eqnarray*}
    X\big(t,Z^0(x)+\tilde{\lambda}_i\,t\big) & = & \int_0^{Z^0(x)+\tilde{\lambda}_i\,t}
\frac{1}{N\big(\tilde{w}^0(\xi)\big)} \, \md\xi
+ \int_0^t \Big(\frac{M}{N}\Big) \big(\tilde{w}\big(\tau,Z^0(x)
+\tilde{\lambda}_i\,t\big)\big)\,\md\tau \\[2mm]
    & = & \int_0^{Z^0(x)} \frac{1}{N\big(\tilde{w}^0(\xi)\big)} \, \md\xi
+ \int_{Z^0(x)}^{Z^0(x)+\tilde{\lambda}_i\,t} \frac{1}{N\big(\tilde{w}^0(\xi)\big)}
\, \md\xi \\[2mm]
    & \; & + \int_0^t \Big(\frac{M}{N}\Big) \big(\tilde{w}\big(\tau,Z^0(x)
+\tilde{\lambda}_i\,t\big)\big)\,\md\tau.
\end{eqnarray*}
A direct computation gives
$$  \int_0^{Z^0(x)} \frac{1}{N\big(\tilde{w}^0(\xi)\big)} \, \md\xi = x.   $$
Together with (\ref{5.1}) at $\bar{w}$, i.e.,
$$  \bar{\lambda}_i = \frac{M}{N}(\bar{w})
- \frac{\tilde{\lambda}_i}{N(\bar{w})}, $$
we obtain (\ref{eq0}). \hfill $\Box$

\begin{lemma}
\label{L3.2}
   It holds
\stepcounter{long}
\begin{equation}
\label{F3.1a}
    \frac{M}{N} = \frac{\tilde{\lambda}_i \lambda_j -\tilde{\lambda}_j \lambda_i}
{\tilde{\lambda}_i- \tilde{\lambda}_j},
\quad \forall \, i, j \in K_n \; \mbox{such that} \; \lambda_j \neq \lambda_i.
\end{equation}
In particular,
\stepcounter{long}
\begin{equation}
\label{F3.1b}
    \frac{M}{N} = \lambda_i \quad \mbox{if} \; \; \tilde{\lambda}_i=0,
\quad \forall \, i \in K_n.
\end{equation}
\end{lemma}

\noindent {\bf Proof.} From (\ref{5.1}), we have
$$ \tilde{\lambda}_i - N \lambda_i = \tilde{\lambda}_j - N \lambda_j,
\quad \forall \, i,j \in K_n.  $$
It follows that
$$  \frac{1}{N} = \frac{\lambda_i -\lambda_j}{\tilde{\lambda}_i-\tilde{\lambda}_j},
\quad \forall \, i, j \in K_n \; \mbox{such that} \; \lambda_j \neq \lambda_i.   $$
Hence, (\ref{5.1}) yields
$$ \frac{M}{N} = \lambda_i - \frac{\tilde{\lambda}_i}{N}
= \frac{\tilde{\lambda}_i \lambda_j -\tilde{\lambda}_j \lambda_i}
{\tilde{\lambda}_i- \tilde{\lambda}_j},
\quad \forall \, i, j \in K_n \; \mbox{such that} \; \lambda_j \neq \lambda_i. $$
This proves the Lemma. \hfill $\Box$

\vspace{3mm}

In the next two lemmas we deal with the limits as $t \rightarrow +\infty$
of the function defined in (\ref{eq0}). For the sake of simplicity, we only
prove the result in the strictly hyperbolic case. In the case of constant
multiplicity of eigenvalues, the result can be obtained in an analogous way
(see Remark \ref{R2.2}).

\begin{lemma}
\label{L3.3}
   Let us denote $u = \displaystyle{\frac{M}{N}}$ and
\stepcounter{long}
\begin{equation}
\label{F3.1c}
  \psi^t_{i1}(x) = \int_0^t \big[u(\tilde{w}(\tau,Z^0(x)+\tilde{\lambda}_i t))
- u(\bar{w})\big] \md \tau, \quad \forall \, x \in \R, \; t > 0.
\end{equation}
Let the assumptions of Proposition \ref{P1.2} hold and the system be strictly
hyperbolic. Let $w^0-\bar{w} \in L^1(\R) \cap L^\infty(\R)$ and
$\dis{\lim_{x \rightarrow \pm \infty} w^0(x) = \bar{w}}$. Then for all $x \in \R$
and $i \in K_n$, we have
$$   \lim_{t \rightarrow +\infty} \psi^t_{i1}(x) = \psi_{i1}(x),  $$
where
\stepcounter{long}
\begin{equation}
\label{F3.1d}
  \hspace{-0.1mm}\psi_{i1}(x) \!=\!\left\{\begin{array}{ll}
  \displaystyle{\sum_{\lambda_j > \lambda_i} \frac{1}{\tilde{\lambda}_j-\tilde{\lambda}_i}
\int_{-\infty}^{+\infty} \!\big(\lambda_i(\bar{w}_1,\cdots,\bar{w}_{j-1},
\tilde{w}^0_j(\xi),\bar{w}_{j+1},\cdots,\bar{w}_n)-\bar{\lambda}_i\big) \md \xi} \\[8mm]
+ \displaystyle{\frac{1}{\tilde{\lambda}_i-\tilde{\lambda}_k} \int_{Z^0(x)}^{+\infty}
\!\big(\lambda_k(\bar{w}_1,\cdots,\bar{w}_{i-1},\tilde{w}^0_i(\xi),\bar{w}_{i+1},
\cdots,\bar{w}_n)-\bar{\lambda}_k\big) \md \xi},
\quad \mbox{if} \; \tilde{\lambda}_i> 0, \\[6mm]
\displaystyle{\int_0^{+\infty} \!\big(\lambda_i(\tilde{w}(\tau,Z^0(x)))
-\lambda_i(\bar{w})\big) \md \tau,
\quad \quad \quad \quad \quad \quad \quad \quad \quad \quad \quad \quad
\quad\; \mbox{if} \; \tilde{\lambda}_i= 0}, \\[4mm]
\displaystyle{\sum_{\lambda_j < \lambda_i} \frac{1}{\tilde{\lambda}_i-\tilde{\lambda}_j}
\int_{-\infty}^{+\infty} \!\big(\lambda_i(\bar{w}_1,\cdots,\bar{w}_{j-1},
\tilde{w}^0_j(\xi),\bar{w}_{j+1},\cdots,\bar{w}_n)-\bar{\lambda}_i\big) \md \xi} \\[8mm]
+ \displaystyle{\frac{1}{\tilde{\lambda}_i-\tilde{\lambda}_k} \int_{Z^0(x)}^{-\infty}
\!\big(\lambda_k(\bar{w}_1,\cdots,\bar{w}_{i-1},\tilde{w}^0_i(\xi),\bar{w}_{i+1},
\cdots,\bar{w}_n)-\bar{\lambda}_k\big) \md \xi},
\quad \mbox{if} \; \tilde{\lambda}_i < 0,
\end{array}\right.\hspace{-1cm}
\end{equation}
with $k \in K_n$ being any index such that $\lambda_k \neq \lambda_i$. Moreover,
the above convergence is uniform on any compact set of $\R$.
\end{lemma}

\noindent {\bf Proof.} Obviously, $\tilde{w}^0-\bar{w} \in L^1(\R)$. Since $Z^0$
is Lipschitzian and strictly increasing on $\R$, we get
$\tilde{w}(\cdot,Z^0(x))-\bar{w} \in L^1(\R)$ for all $x \in \R$. Then the
Taylor formula implies that
$\lambda_i(\tilde{w}(\cdot,Z^0(x)))-\lambda_i(\bar{w}) \in L^1(\R)$.
Hence, the result for $\tilde{\lambda}_i= 0$ follows from (\ref{F3.1b}).

Now we prove the result for $\tilde{\lambda}_i> 0$ and we omit the proof of the
result for $\tilde{\lambda}_i< 0$ since it is similar. First, the Taylor
formula gives
$$  u\big(\tilde{w}(\tau,Z^0(x)+\tilde{\lambda}_i t)\big)- u(\bar{w})
= \sum_{j=1}^n u^t_{ij}(\tau,x) \big(\tilde{w}_j(\tau,Z^0(x)
+ \tilde{\lambda}_i t)-\bar{w}_j\big),  $$
where
$$  u^t_{ij}(\tau,x) = \int_0^1 \frac{\pa u}{\pa w_j}\big(\bar{w}
+\alpha(\tilde{w}(\tau,Z^0(x)+\tilde{\lambda}_i t)- \bar{w})\big) \md \alpha.  $$
If $\tilde{\lambda}_j=0$ for some $j \in K_n$, from (\ref{F3.1b}) we have
$u = \lambda_j$. Since the system is linearly degenerate, $\lambda_j$ is independent
of $w_j$, so that $\dis \frac{\pa u}{\pa w_j} = 0$. Hence, from (\ref{tdw}) we may write
$$ u\big(\tilde{w}(\tau,Z^0(x)+\tilde{\lambda}_i t)\big) - u(\bar{w})
 = \sum_{\tilde{\lambda}_j \neq 0} u^t_{ij}(\tau,x)
\big(\tilde{w}^0_j(Z^0(x)+\tilde{\lambda}_i t-\tilde{\lambda}_j \tau)-\bar{w}_j\big).
$$
For each $\tilde{\lambda}_j \neq 0$, we may make a change of variable
$\xi = Z^0(x) + \tilde{\lambda}_i t-\tilde{\lambda}_j \tau$ in the integral of
the above functions. It follows that
$$   \psi^t_{i1}(x) = - \sum_{\tilde{\lambda}_j \neq 0} \frac{1}{\tilde{\lambda}_j}
\int_{Z^0(x)+\tilde{\lambda}_i\,t}^{Z^0(x)+(\tilde{\lambda}_i-\tilde{\lambda}_j)t}
u^t_{ij}\big(\frac{Z^0(x)+\tilde{\lambda}_i t-\xi}{\tilde{\lambda}_j},x\big)
\big(\tilde{w}^0_j(\xi)-\bar{w}_j\big) \md \xi.     $$

Note that the function $u^t_{ij}$ is uniformly bounded, i.e., there is a
constant $C>0$ such that
$$  |u^t_{ij}(\tau,x)| \leq C, \quad \forall \, \tau \in [0,t], \; x \in \R, \; t> 0.  $$
Since $w^0-\bar{w} \in L^1(\R)$, we obtain
$$ |u^t_{ij}(\tilde{w}^0_j-\bar{w}_j)| \leq C |\tilde{w}^0_j-\bar{w}_j| \in L^1(\R).  $$
On the other hand, for all $k \in K_n$, from
$\xi = Z^0(x) + \tilde{\lambda}_i t-\tilde{\lambda}_j \tau$ we have
$$   Z^0(x) + \tilde{\lambda}_i t-\tilde{\lambda}_k \tau
= \frac{\tilde{\lambda}_k}{\tilde{\lambda}_j}\xi
+ \Big(1-\frac{\tilde{\lambda}_k}{\tilde{\lambda}_j}\Big)
(Z^0(x) + \tilde{\lambda}_i t).    $$
Since $Z^0$ is a Lipschitz function, for all $x \in \R$ and $\xi \in \R$, we have
$$   \lim_{t \rightarrow +\infty} \big(Z^0(x) + \tilde{\lambda}_i t
-\tilde{\lambda}_k \tau\big) = \infty,
\quad \forall \, k \in K_n \; \mbox{such that} \; \lambda_k \neq \lambda_j.   $$
It follows from (\ref{tdw}) and $\dis{\lim_{x \rightarrow \pm \infty} w^0(x) = \bar{w}}$
that
$$   \lim_{t \rightarrow +\infty} \big(\tilde{w}^0_k(Z^0(x)
+ \tilde{\lambda}_i t-\tilde{\lambda}_k \tau)\big)
= \bar{w}_k, \quad \forall \, k \in K_n \; \mbox{such that}
\; \lambda_k \neq \lambda_j.   $$
Hence,
$$    \lim_{t \rightarrow +\infty} u^t_{ij}\big(\frac{Z^0(x)+\tilde{\lambda}_i t
-\xi}{\tilde{\lambda}_j},x\big)
= \int_0^1 \frac{\pa u}{\pa w_j}\big(\bar{w}_1,\cdots,\bar{w}_{j-1},\bar{w}_j
+\alpha(\tilde{w}^0_j(\xi) -\bar{w}_j),\bar{w}_{j+1},\cdots,\bar{w}_n\big)\, \md \alpha.  $$
Or
\begin{eqnarray}
    & \, & \int_0^1 \frac{\pa u}{\pa w_j}\big(\bar{w}_1,\cdots,\bar{w}_{j-1},\bar{w}_j
+\alpha(\tilde{w}^0_j(\xi) -\bar{w}_j),\bar{w}_{j+1},\cdots,\bar{w}_n\big)\, \md \alpha \cdot
\big(\tilde{w}^0_j(\xi)-\bar{w}_j\big) \nonumber \\[1mm]
    & = & u(\bar{w}_1,\cdots,\bar{w}_{j-1},\tilde{w}^0_j(\xi) ,\bar{w}_{j+1},\cdots,\bar{w}_n) - u(\bar{w}). \nonumber
\end{eqnarray}

For $\tilde{\lambda}_j < \tilde{\lambda}_i$, we have
$Z^0(x) + (\tilde{\lambda}_i-\tilde{\lambda}_j)t \rightarrow +\infty$. Since
$\tilde{\lambda}_i > 0$, we also have $Z^0(x) + \tilde{\lambda}_i t \rightarrow +\infty$.
Therefore, for all $x \in \R$, the Lebesgue dominated convergence theorem yields
$$  \lim_{t \rightarrow +\infty}
\int_{Z^0(x)+\tilde{\lambda}_it}^{Z^0(x)+(\tilde{\lambda}_i-\tilde{\lambda}_j)t}
\,u^t_{ij}\big(\frac{Z^0(x)+\tilde{\lambda}_i t-\xi}{\tilde{\lambda}_j},x\big)
\big(\tilde{w}^0_j(\xi)-\bar{w}_j\big)\, \md \xi = 0.   $$
Similarly,
\begin{eqnarray}
    & \, & \lim_{t \rightarrow +\infty}
    \int_{Z^0(x)+\tilde{\lambda}_it}^{Z^0(x)+(\tilde{\lambda}_i-\tilde{\lambda}_j)t}
u^t_{ij}\big(\frac{Z^0(x)+\tilde{\lambda}_i t-\xi}{\tilde{\lambda}_j},x\big)
\big(\tilde{w}^0_j(\xi)-\bar{w}_j\big) \md \xi \nonumber \\[1mm]
    & = & \int_{+\infty}^{-\infty}
\big[u(\bar{w}_1,\cdots,\bar{w}_{j-1},\tilde{w}^0_j(\xi) ,\bar{w}_{j+1},\cdots,\bar{w}_n)
- u(\bar{w})\big] \md \xi,
\quad \mbox{if} \quad \tilde{\lambda}_j > \tilde{\lambda}_i, \nonumber
\end{eqnarray}
and
\begin{eqnarray}
    & \, & \lim_{t \rightarrow +\infty}
    \int_{Z^0(x)+\tilde{\lambda}_it}^{Z^0(x)+(\tilde{\lambda}_i-\tilde{\lambda}_j)t}
u^t_{ij}\big(\frac{Z^0(x)+\tilde{\lambda}_i t-\xi}{\tilde{\lambda}_j},x\big)
\big(\tilde{w}^0_j(\xi)-\bar{w}_j\big) \md \xi \nonumber \\[1mm]
    & = & \int_{+\infty}^{Z^0(x)}
\big[u(\bar{w}_1,\cdots,\bar{w}_{i-1},\tilde{w}^0_i(\xi) ,\bar{w}_{i+1},\cdots,\bar{w}_n)
- u(\bar{w})\big] \md \xi,
\quad \mbox{if} \quad \tilde{\lambda}_j = \tilde{\lambda}_i. \nonumber
\end{eqnarray}
This implies that
$$   \lim_{t \rightarrow +\infty} \psi^t_{i1}(x) = \psi_{i1}(x)  $$
with
\begin{eqnarray}
    \psi_{i1}(x) & = &
  \displaystyle{\sum_{\tilde{\lambda}_j > \tilde{\lambda}_i} \frac{1}{\tilde{\lambda}_j}
\int_{-\infty}^{+\infty} \!\big[u(\bar{w}_1,\cdots,\bar{w}_{j-1},
\tilde{w}^0_j(\xi),\bar{w}_{j+1},\cdots,\bar{w}_n)-u(\bar{w})\big] \md \xi} \nonumber \\[1mm]
& \; & +\; \displaystyle{\frac{1}{\tilde{\lambda}_i} \int_{Z^0(x)}^{+\infty}
\!\big[u(\bar{w}_1,\cdots,\bar{w}_{i-1},\tilde{w}^0_i(\xi),\bar{w}_{i+1},
\cdots,\bar{w}_n)-u(\bar{w})\big] \md \xi}. \nonumber
\end{eqnarray}

Finally, from (\ref{5.1}) and $N> 0$, we deduce that $\lambda_i - \lambda_j$
and $\tilde{\lambda}_i -\tilde{\lambda}_j$ have the same sign. Using (\ref{F3.1a})
and the fact that $\lambda_j$ is independent of $w_j$ for all $j \in K_n$, we
obtain (\ref{F3.1d}) for $\tilde{\lambda}_i >0$. This proves the pointwise
convergence of $\psi^t_{i1}$ on $\R$ for all $i \in K_n$. When $x$ stays in
a compact set of $\R$, $Z^0(x)$ is bounded independent of $x$. This shows that
the convergence of $\psi^t_{i1}$ is uniform on any compact set of $\R$.
\hfill $\Box$

\begin{lemma}
\label{L3.4}
    Under the assumptions of Lemma \ref{L3.3}. For each $i \in K_n$ and
all $x \in \R$, we have
\stepcounter{long}
\begin{equation}
\label{eq3}
    \lim_{t \rightarrow +\infty} \big(X(t,Z^0(x)+\tilde{\lambda}_i\,t)
-\bar{\lambda}_i\,t\big) = \psi_i(x),
\end{equation}
where
\stepcounter{long}
\begin{equation}
\label{eq4}
    \psi_i(x)=x + \psi_{i1}(x) + \psi_{i2}(x),
\end{equation}
with
\stepcounter{long}
\begin{equation}
\label{eq4a}
    \psi_{i2}(x) = \left\{\begin{array}{ll}
\displaystyle{\int_{Z^0(x)}^{+\infty}
\Big(\frac{1}{N\big(\tilde{w}^0(\xi)\big)} - \frac{1}{N(\bar{w})}\Big)\, \md\xi},
\; \; & \mbox{if} \quad \tilde{\lambda}_i> 0, \\[5mm]
\quad \quad 0, \; \; & \mbox{if} \quad \tilde{\lambda}_i= 0, \\[3mm]
\displaystyle{\int_{Z^0(x)}^{-\infty}
\Big(\frac{1}{N\big(\tilde{w}^0(\xi)\big)} - \frac{1}{N(\bar{w})}\Big)\, \md\xi},
\; \; & \mbox{if} \quad \tilde{\lambda}_i< 0.
\end{array}\right.
\end{equation}
Moreover, the convergence (\ref{eq3}) is uniform on any compact set of $\R$.
\end{lemma}

\noindent {\bf Proof.} From $w^0 -\bar{w} \in L^1(\R) \cap L^\infty(\R)$, we have
$$   \frac{1}{N(\tilde{w}^0)} -\frac{1}{N(\bar{w})} \in L^1(\R),   $$
so that, for all $x \in \R$, when $\tilde{\lambda}_i > 0$,
$$  \lim_{t \rightarrow +\infty} \int_{Z^0(x)}^{Z^0(x) + \tilde{\lambda}_i t}
\left(\frac{1}{N(\tilde{w}^0)} -\frac{1}{N(\bar{w})}\right)\!\md \xi
= \int_{Z^0(x)}^{+\infty} \left(\frac{1}{N(\tilde{w}^0)}
-\frac{1}{N(\bar{w})}\right)\!\md \xi.   $$
Thus, the result for $\tilde{\lambda}_i > 0$ follows from Lemmas \ref{L3.1}
and \ref{L3.3}. The cases $\tilde{\lambda}_i < 0$ and $\tilde{\lambda}_i = 0$
can be treated in a same way. \hfill $\Box$

\vspace{3mm}

The next two results are concerned with regularities and monotonicities of
functions $\psi_i, \; i \in K_n$. To this end, we distinguish again the cases
$\tilde{\lambda}_i \neq 0$ and $\tilde{\lambda}_i = 0$. In both cases, we need
smallness conditions on the initial data to only assure the invertibility of
$\psi_i$.

\begin{lemma}
\label{L3.5}
   Let the assumptions of Lemma \ref{L3.3} hold. For each $i \in K_n$
with $\tilde{\lambda}_i \neq 0$, $\psi_i$ is a Lipschitz function on $\R$.
Moreover, there are constants $\delta_1> 0$ and $C_{\delta_1}>0$ such that if
$\|w^0-\bar{w}\|_{L^\infty(\R)} < \delta_1$, then $\psi'_i(x) \geq C_{\delta_1}$
for almost all $x \in \R$.
\end{lemma}

\noindent {\bf Proof.} We only prove the result for $i \in K_n$ with
$\tilde{\lambda}_i >0$. The case for $\tilde{\lambda}_i <0$ can be treated
in a same way.

From (\ref{F3.1d}) and (\ref{eq4})-(\ref{eq4a}), in the strictly hyperbolic
case, we have
\stepcounter{long}
\begin{eqnarray}
\label{eq4b}
  \psi_i(x) & = & x + \int_{Z^0(x)}^{+\infty} \Big(\frac{1}{N(\tilde{w}^0(\xi))}
-\frac{1}{N(\bar{w})}\Big) \md \xi \\[2mm]
& \, & + \frac{1}{\tilde{\lambda}_i-\tilde{\lambda}_k}
\int_{Z^0(x)}^{+\infty} \big(\lambda_k(\bar{w}_1,\cdots,\bar{w}_{i-1},\tilde{w}^0_i(\xi),
\bar{w}_{i+1},\cdots,\bar{w}_n) -\bar{\lambda}_k\big)\md \xi + \bar{\psi}_i, \nonumber
\end{eqnarray}
for any $k \in K_n$ such that $\tilde{\lambda}_k \neq \tilde{\lambda}_i$, where
$\bar{\psi}_i$ is a constant. Since $Z^0$ is a Lipschitz function on $\R$ and
the function in the integral of (\ref{eq4b}) is in $L^1(\R)\cap L^{\infty}(\R)$,
by composition we deduce that $\psi_i$ is a Lipschitz function on $\R$.

Hence, $\psi_i$ is differentiable for almost all $x \in \R$. Let us denote
$$  f_i(w) = \frac{1}{N(w)} + \frac{\lambda_k(\bar{w}_1,\cdots,\bar{w}_{i-1},w_i,
\bar{w}_{i+1},\cdots,\bar{w}_n)}{\tilde{\lambda}_i-\tilde{\lambda}_k},  $$
which is a smooth function. Noting $\tilde{w}^0(Z^0) = w^0$, a straightforward
computation yields
$$  \psi'_i(x) = 1 - N(w^0(x)) \big(f_i(w^0(x)) - f_i(\bar{w})\big).    $$
Thus, when $\|w^0-\bar{w}\|_{L^\infty(\R)} < \delta_1$, which is sufficiently small,
we obtain $\psi'_i(x) \geq C_{\delta_1} > 0$ for almost all $x \in \R$. This proves
the result. \hfill $\Box$

\begin{lemma}
\label{L3.6}
    Let the assumptions of Lemma \ref{L3.3} hold. For each $i \in K_n$ with
$\tilde{\lambda}_i = 0$, $\psi_i$ is a uniform continuous function on $\R$.
If $w^0 -\bar{w} \in W^{1,1}(\R)$, then $\psi_i$ is a Lipschitz function on $\R$.
Moreover, there are constants $\delta_2> 0$ and $C_{\delta_2} > 0$ such that if
$\|(w^0)'\|_{L^1(\R)} < \delta_2$, then $\psi'_i(x) \geq C_{\delta_2}$ for almost all
$x \in \R$.
\end{lemma}

\noindent {\bf Proof.} Replacing $w^0-\bar{w}$ by $w^0$, we may suppose
$\bar{w} = 0$. Since $\tilde{\lambda}_i = 0$, from (\ref{F3.1d}) and
(\ref{eq4})-(\ref{eq4a}), we have
$$  \psi_i(x) = x + \int_0^{+\infty} \!\big(\lambda_i(\tilde{w}(\tau,Z^0(x)))
-\lambda_i(0)\big) \md \tau.  $$
Let $x$, $y \in \R$. Then
\begin{eqnarray}
    \psi_i(x) - \psi_i(y) & = & (x-y)+\int_0^{+\infty} \!\Big(\lambda_i\big(\tilde{w}^0_1(Z^0(x)
-\tilde{\lambda}_1 \tau),\cdots,\tilde{w}^0_n(Z^0(x)
-\tilde{\lambda}_n\tau )\big) \nonumber \\[2mm]
& \, & -\lambda_i\big(\tilde{w}^0_1(Z^0(y)-\tilde{\lambda}_1 \tau),
\cdots,\tilde{w}^0_n(Z^0(y) -\tilde{\lambda}_n \tau)\big)\Big)\, \md \tau. \nonumber
\end{eqnarray}
Let us denote
$$   \Lambda= \max\Big\{\Big|\frac{\pa \lambda_i(w)}{\pa w_j}\Big|;
\; w \in L^\infty(\R), \; \forall \, i, \, j \in K_n\Big\}.             $$
By the Taylor formula, we get
$$  \big|\psi_i(x) - \psi_i(y)\big| \leq |x-y|+ \Lambda \sum_{\tilde{\lambda}_j \neq \tilde{\lambda}_i}
\int_0^{+\infty} \big|\tilde{w}^0_j(Z^0(x)-\tilde{\lambda}_j \tau)
- \tilde{w}^0_j(Z^0(y)-\tilde{\lambda}_j \tau)\big| \md \tau.   $$
Noting that $Z^0$ is a Lipschitz function and $w^0 \in L^1(\R)$, from the uniform
continuity of function $x \longmapsto w^0(x+\cdot)$ in $L^1(\R)$ (see, for instance,
\cite{Br83}) and by composition, we conclude that $\psi_i$ is uniform continuous
on $\R$.

Suppose now $w^0 \in W^{1,1}(\R)$. It is easy to check that $\psi_i$ is a Lipschitz
function on $\R$ and we have
$$      \psi_i'(x) = 1 + \sum_{\tilde{\lambda}_j \neq 0}
\int_0^{+\infty}\, \frac{\pa \lambda_i}{\pa w_j}\cdot \frac{\pa \tilde{w}_j^0(Z^0(x)
-\tilde{\lambda}_j\,\tau)}{\pa x} \,\md \tau.     $$
By changes of variables, we get
$$    \psi_i'(x) = 1 - \sum_{\tilde{\lambda}_j > 0} \frac{1}{\tilde{\lambda}_j}
\int_{Z^0(x)}^{-\infty} \, \frac{\pa \lambda_i}{\pa w_j}
\big(\tilde{w}_j^0(\xi)\big)' N(w^0(x)) \,\md \xi
- \sum_{\tilde{\lambda}_j < 0} \frac{1}{\tilde{\lambda}_j}
\int_{Z^0(x)}^{+\infty} \, \frac{\pa \lambda_i}{\pa w_j}
\big(\tilde{w}_j^0(\xi)\big)' N(w^0(x)) \,\md \xi.    $$
Therefore, there is a constant $C_1>0$ such that
$$  \psi_i'(x) \geq 1 - C\|(w^0)'\|_{L^1(\R)}, \quad \forall \, x \in \R.   $$
Thus, when $\|(w^0)'\|_{L^1(\R)} < \delta_2$, which is sufficiently small, we get
$\psi'_i(x) \geq C_{\delta_2} > 0$ for almost all $x \in \R$. This completes the proof of
Lemma \ref{L3.6}. \hfill $\Box$

\vspace{3mm}

To state and prove the main result of this section, we still need a
lemma, which was established by the authors in \cite{PY10}.

\begin{lemma}
\label{L3.7}
    Let $f - \bar{f} \in L^1(\R)$ with $\bar{f}$ being a real constant. Let
$(e_\nu)_{\nu > 0}$ be a sequence of Lipschitz functions defined on $\R$.
Suppose that there is a constant $\gamma> 0$ such that $e'_\nu(x) \geq \gamma$,
$\forall \, \nu> 0$, a.e. $x \in \R$ and $e_\nu(x) \longrightarrow x$ as
$\nu \rightarrow 0$, pointwisely on $\R$ and uniformly on any compact set of $\R$.
Then
$$I^\nu_f = \int_{\R}\, \big|f(x) - f(e_\nu(x))\big|\, \md x \longrightarrow 0,
\quad \text{as} \quad \nu \rightarrow 0.
$$
\end{lemma}

\begin{theorem}
\label{T3.1}
      Let the assumptions of Lemma \ref{L3.3} hold. Furthermore, assume that for
each $i \in K_n$, $\psi_i$ is Lipschitzian on $\R$ satisfying $\psi'_i(x) \geq C_0$
for some constant $C_0>0$. Let $\phi_i$ be the inverse function of $\psi_i$. Then
we have
\stepcounter{long}
\begin{equation}
\label{F3.10a}
  \dis \lim_{t \rightarrow +\infty} \|w_i(t,\cdot)-w_i^0(\phi_i(\cdot
-\bar{\lambda}_i\,t))\|_{L^1(\R)}=0.
\end{equation}
\end{theorem}

\noindent {\bf Proof.} Let $x=X^0(Z(t,\eta)-\tilde{\lambda}_i\,t)$.
Then $\eta=X(t,Z^0(x)+\tilde{\lambda}_i\,t)$. From
(\ref{Z0})-(\ref{X1}), we have
$$   \dis \md \eta=J^0(t,x)\,\md x,\quad \quad \text{with}\quad
J^0(t,x)=\frac{N(w^0(x))}{N(\tilde{w}(t,Z^0(x)+\tilde{\lambda}_i\,t))}.  $$
Obviously, there are positive constants denoted by $C_1$ and $C_2$ such that
$C_1 \leq J^0(t,x) \leq C_2$ for almost all $(t,x)\in \R_+\times\R$. We deduce
from the explicit expression of the entropy solutions that
\begin{align}
\dis\;&\;\|w_i(t,\cdot)-w_i^0(\phi_i(\cdot-\bar{\lambda}_i\,t))\|_{L^1(\R)}\nonumber\\
= \;&\; \int_{\R}\,|w_i^0(X^0(Z(t,\eta)-\tilde{\lambda}_i\,t))
-w_i^0(\phi_i(\eta-\bar{\lambda}_i\,t))|\,\md \eta \nonumber \\
= \;&\; \int_{\R}\,|w_i^0(x)-w_i^0(\phi_i(X(t,Z^0(x)+\tilde{\lambda}_i\,t)
-\bar{\lambda}_i\,t))|J^0(t,x)\,\md x. \nonumber
\end{align}

Let
$$e_t(x)=\phi_i(X(t,Z^0(x)+\tilde{\lambda}_i\,t)-\bar{\lambda}_i\,t).$$
It is easy to see that function $w_i^0-w_i^0(e_t)$ is uniformly
bounded in $L^1(\R)$ with respect to $t$ and
$$  e_t'(x)= J^0(t,x)\,\phi_i'(X(t,Z^0(x)
+\tilde{\lambda}_i\,t)-\bar{\lambda}_i\,t).  $$
From the discussion above, there is a positive constant, denoted by
$\gamma$, such that $e_t'(x) \geq \gamma$ for all $(t,x)\in
\R_+\times \R$. Since $\phi_i$ is the inverse function of $\psi_i$,
Lemma \ref{L3.4} yields
$$   \forall\,x\in \R,\quad \quad \lim_{t \rightarrow +\infty} e_t(x)=x,  $$
pointwisely on $\R$ and uniformly on any compact set of $\R$. Thus, $e_t$ fulfills the conditions of Lemma \ref{L3.7} with
$\Omega_1=\R$ and $\dis \nu=\frac{1}{t}$. We conclude that
$$  \lim_{t \rightarrow +\infty} \int_{\R}\,|w_i^0(x)-w_i^0(\phi_i(X(t,Z^0(x)
+\tilde{\lambda}_i\,t) -\bar{\lambda}_i\,t))|J^0(t,x)\,\md x = 0,   $$
which completes the proof of Theorem \ref{T3.1}. \hfill $\Box$

\begin{remark}
\label{R3.1}
    In Theorem \ref{T3.1}, the assumption $\psi'_i(x) \geq C_0$ for all
$x \in \R$ is essential. It implies that $\psi_i$ is a strictly increasing
function on $\R$. Sufficient conditions to ensure this assumption are
given in Lemmas \ref{L3.5} and \ref{L3.6}, from which we see that the case
$\tilde{\lambda_i} = 0$ requires more regularity on the initial data. In
particular, if $\tilde{\lambda}_i \neq 0$ for some $i \in K_n$, we only require
that $w^0-\bar{w} \in L^1(\R) \cap L^\infty(\R)$ and the $L^\infty$ norm of
$w^0-\bar{w}$ is small to yield (\ref{F3.10a}).
\end{remark}

%$4%%%%%%%%%%%%%%%%%%%%%%%%%%%%%%%%%%%%%%%%%%%%%%%%%%%%%%%%%%%%%%%%%%%%%%

\section{$L^1$ stability of entropy solutions}
\newcounter{exte}
\renewcommand{\theequation}{\thesection.\theexte}

This section is concerned with the $L^1$ stability of entropy solutions to the
Cauchy problem (\ref{diag})-(\ref{diag-init}) in the case $N_i =N$ for all
$i \in K_n$. To this end, suppose $w^{1,0}$, $w^{2,0} \in L^1(\R)\cap L^{\infty}(\R)$
and let $w^1$ and $w^2$ be the entropy solutions with initial data $w^{1,0}$
and $w^{2,0}$, respectively. For simplicity, we define
$$    R_0=\|w^{2,0}-w^{1,0}\|_{L^1(\R)}\quad\text{and}
\quad R_t=\|w^{2}(t,\cdot)-w^{1}(t,\cdot)\|_{L^1(\R)}.      $$
Here the $L^1$ stability of the entropy solutions to the Cauchy problem
(\ref{diag})-(\ref{diag-init}) means that for all $\eps>0$,
there exists $\delta > 0$ depending only on $\eps$ such that
$$  R_0 \leq \delta \Longrightarrow R_t \leq \eps, \quad \forall \; t > 0,  $$
or equivalently,
\stepcounter{exte}
\begin{equation}
\label{rtr0}
    \lim_{R_0 \rightarrow 0} R_t = 0, \quad \forall \; t > 0.
\end{equation}
Together with the uniqueness of entropy solutions, stability
condition (\ref{rtr0}) implies that application $R_0 \longmapsto
R_t$ is continuous at $R_0 = 0$ for all $t \geq 0$. In addition, the
strong $L^1$ stability means that the following stability inequality
holds for some constant $C>0$ independent of $t$:
\stepcounter{exte}
\begin{equation}
\label{rtr1}
    R_t \leq C R_0, \quad \forall \; t > 0.
\end{equation}
The study of the strong $L^1$ stability needs usually to work with
small BV initial data, see \cite{Br00} for a general system of
conservation laws and \cite{Ch92} for system (\ref{diag}).

For each $j =1, 2$, since $w^{j,0}\in L^1(\R)\cap L^{\infty}(\R)$,
from Proposition \ref{P1.2}, the Cauchy problem
(\ref{diag})-(\ref{diag-init}) with the initial data $w^{j,0}$
admits a unique entropy solution, given by
$$    w^j(t,x)=\Big(w_1^{j,0}\big(X^{j,0}(Z^j(t,x)-\tilde{\lambda}_1\,t)\big),
\cdots,w_n^{j,0}\big(X^{j,0}(Z^j(t,x)-\tilde{\lambda}_n\,t)\big)\Big)^{\top}, $$
where $X^{j,0}$ is the inverse function of $Z^{j,0}$ and
$Z^j(t,\cdot)$ is the inverse function of $X^j(t,\cdot)$, in which
$$ Z^{j,0}(x)=\int_0^x\,N\big(w^{j,0}(\xi)\big)\,\md \xi $$
and
$$       \md X^j(t,z)=\frac{1}{N(\tilde{w}^j(t,z))}\,\md z
+\frac{M(\tilde{w}^j(t,z))}{N(\tilde{w}^j(t,z))}\,\md t,  $$
with
\stepcounter{exte}
\begin{equation}
\label{tlwj}
\tilde{w}^j(t,z)=\Big(w_1^{j,0}\big(X^{j,0}(z-\tilde{\lambda}_1\,t)\big)
,\cdots,w_n^{j,0}\big(X^{j,0}(z-\tilde{\lambda}_n\,t)\big)\Big)^{\top}.
\end{equation}

The main result of this section is as follows.

\begin{theorem}
\label{T4.1}
Under the assumptions of Proposition \ref{P1.2}, the unique global entropy
solution to the Cauchy problem (\ref{diag})-(\ref{diag-init}) is
$L^1$ stable. If furthermore, $w^{j,0}\in W^{1,1}(\R)$ for $j=1,2$,
then, for each $i\in K_n$, we have
\stepcounter{exte}
\begin{equation}
\label{stability'}
  \|w_i^2(t,\cdot)-w_i^1(t,\cdot)\|_{L^1(\R)} \leq C \; R_0, \quad \forall \,t \geq 0.
\end{equation}
\end{theorem}

The proof of Theorem \ref{T4.1} is based on a series of Lemmas below.

\begin{lemma}
\label{L5.2}
       Let the assumptions of Proposition \ref{P1.2} hold. There exists a
positive constant $C$, such that for all $x$, $z \in \R$, we have
$$  \left|Z^{2,0}(x)-Z^{1,0}(x)\right|\leq C R_0, \quad
 \left|X^{2,0}(z)-X^{1,0}(z)\right|\leq C R_0  $$
and
\stepcounter{exte}
\begin{equation}
\label{5.0}
    \left|X^{1,0}(Z^{2,0}(x))-x\right| \leq C R_0.
\end{equation}
\end{lemma}

\noindent {\bf Proof.} The proof is omitted here since it is similar to that
of Lemma 4.1 in \cite{PY10}. \hfill $\Box$

\vspace{3mm}

The next lemma was proved by the authors in \cite{PY10}.

\begin{lemma}
\label{L4.2}
    Let $f \in L^1_{loc}(\R)$, $f' \in L^1(\R)$  and $h_t$ be a Lipschitz function
on $\R$ such that
$$  \text{a.e.} \; (t,x) \in \R_+ \times \R, \quad h'_t(x) \geq \gamma_1,
\quad \big|h_t(x) - x\big| \leq \gamma_2(t),  $$ for a constant
$\gamma_1 > 0$ and a function $\gamma_2(t) \geq 0$. Then $f-f(h_t)
\in L^1(\R)$ and we have
$$    \int_{\R}\, \big|f(x) -f(h_t(x))\big|\, \md x \leq
\frac{\gamma_2(t)}{\gamma'_1}\; \|f'\|_{L^1(\R)},   $$
where $\gamma'_1 = \min\{1,\gamma_1\}$.
\end{lemma}

\begin{lemma}
\label{L5.3}
    Under the assumptions of Proposition \ref{P1.2}, we have
\stepcounter{exte}
\begin{equation}
\label{5a}
     \left|X^2-X^1\right|\longrightarrow 0 \quad \text{and}
\quad \left|Z^2-Z^1\right|\longrightarrow 0, \quad
\text{as}\quad R_0\rightarrow 0,
\end{equation}
uniformly on $\R^+\times \R$.
If furthermore, $w^{j,0}\in W^{1,1}(\R)$ for $j=1,2$, then, for all
$t \geq 0$ and all $x, \, z \in \R$, we have \stepcounter{exte}
\begin{equation}
\label{conver'}
    \left|X^2(t,z)-X^1(t,z)\right| \leq C\,R_0, \quad
\left|Z^2(t,x)-Z^1(t,x)\right| \leq CR_0.
\end{equation}
\end{lemma}

\noindent {\bf Proof.} For $x \in \R$, let $z_2=Z^2(t,x)$. Then
$x=X^2(t,z_2)$. Since $Z^1(t,\cdot)$ is Lipschitzian on $\R$, we have
\begin{align}
\nonumber
       \left |Z^2(t,x)-Z^1(t,x)\right |=&\left |Z^2(t,X^2(t,z_2))
-Z^1(t,X^2(t,z_2))\right |\\\nonumber
  = &\left |Z^1(t,X^1(t,z_2))-Z^1(t,X^2(t,z_2))\right |\\\nonumber
\leq & \,C \left |X^2(t,z_2)-X^1(t,z_2)\right |.
\end{align}
Hence, the second estimate of (\ref{5a}) immediately follows
from the first one. It remains to prove the first estimate of (\ref{5a}).

From (\ref{X2}), for any $i\in K_n$ we have
$$   X^j(t,z)=\int_0^{z}\,\frac{\md \xi}{N(w^{j,0}(X^{j,0}(\xi)))}
+\int_0^t\lambda_i(\tilde{w}^j(\tau,z))\,\md \tau
-\int_0^t\frac{\tilde{\lambda}_i\,\md \tau}{N(\tilde{w}^j(\tau,z))}, $$
where $\tilde{w}^j(t,z) \; (j=1,2)$ is defined by (\ref{tlwj}). Therefore
\begin{align}
\nonumber
\left |X^2(t,z)-X^1(t,z)\right |\le &\int_{\R}\left|\frac{1}{N(w^{2,0}(X^{2,0}(\xi)))}
-\frac{1}{N(w^{1,0}(X^{1,0}(\xi)))}\right|\,\md \xi\\\nonumber
 &\,+\int_{\R^+}\,\left|\lambda_i(\tilde{w}^2(\tau,z))
-\lambda_i(\tilde{w}^1(\tau,z))\right|\,\md \tau\\\nonumber
 &\,+\tilde{\lambda}_i\int_{\R^+}\left|\frac{1}{N(\tilde{w}^2(\tau,z))}
-\frac{1}{N(\tilde{w}^1(\tau,z))}\right|\md \tau\\
\stackrel{def}{=} &\; {\cal F} + {\cal G}_i + {\cal H}_i. \nonumber
\end{align}
Since $N>0$, we have
\stepcounter{exte}
\begin{align}
\nonumber
 {\cal F}=&\int_{\R}\frac{\left|N(w^{2,0}(X^{2,0}(\xi)))
-N(w^{1,0}(X^{1,0}(\xi)))\right|}{N(w^{1,0}(X^{1,0}(\xi)))\cdot N(w^{2,0}(X^{2,0}(\xi)))}
 \,\md \xi\\\nonumber
 \leq &\,C\, \int_{\R}\left|w^{2,0}(X^{2,0}(\xi))- w^{1,0}(X^{1,0}(\xi))\right|\, \md \xi
 \\\nonumber
 \leq &\, C\, \int_{\R}\left|w^{2,0}(X^{2,0}(\xi))- w^{1,0}(X^{2,0}(\xi))\right|\, \md
 \xi\\\label{i1}
 &\,+\,C\,\int_{\R}\left|w^{1,0}(x)- w^{1,0}(X^{1,0}(Z^{2,0}(x)))\right|\, \md
 x.
 \end{align}
It follows from (\ref{5.0}) that
$$ X^{1,0}(Z^{2,0}(x))\longrightarrow x, \quad \text{as}\;\; R_0\rightarrow 0,$$
uniformly on $\R$. Moreover,
\stepcounter{exte}
\begin{equation}
\label{5.2}
 \forall\; x \in \R,\quad \Big(X^{1,0}(Z^{2,0}(x))\Big)'\,
=\,\dfrac{N(w^{2,0}(x))}{N(w^{1,0}(X^{1,0}(Z^{2,0}(x))))}\geq
\gamma,
\end{equation}
for some constant $\gamma > 0$. We conclude
from Lemma \ref{L3.7} with $\Omega_1=\R$ that
$$ \int_{\R}\left|w^{1,0}(x)- w^{1,0}(X^{1,0}(Z^{2,0}(x)))\right|\, \md
 x\longrightarrow 0,\quad \text{as}\; \;R_0\rightarrow 0,  $$
which together with (\ref{i1}) implies that ${\cal F}\longrightarrow 0$ as
$R_0\rightarrow 0$.

Next, we consider ${\cal G}_i$ and ${\cal H}_i$. If
$\tilde{\lambda}_k \neq 0$ for all $k \in K_n$, then for any fixed $i \in K_n$,
$\mu_i=\lambda_i$ or $\mu_i=\tilde{\lambda}_i/N$, we have
\begin{align}
      \displaystyle & \int_{\R^+}\,\left|\mu_i(\tilde{w}^2(\tau,z))
-\mu_i(\tilde{w}^1(\tau,z))\right|\,\md \tau \nonumber \\
\displaystyle \leq &\; C\;\sum_{k=1}^n\;\int_{\R^+}
\;\left|w_k^{2,0}(X^{2,0}(z-\tilde{\lambda}_k\,\tau))
-w_k^{1,0}(X^{1,0}(z-\tilde{\lambda}_k\,\tau))\right|\;\md \tau, \nonumber
\end{align}
which can be treated as for ${\cal F}$ due to changes of variables
$\tau \longmapsto z-\tilde{\lambda}_k\,\tau$. This shows
${\cal G}_i \longrightarrow 0$ and ${\cal H}_i \longrightarrow 0$ as
$R_0 \rightarrow 0$. If $\tilde{\lambda}_i=0$ for
some $i \in K_n$, then ${\cal H}_i=0$. Since $\lambda_i$ is linearly degenerate,
$\lambda_i$ is independent of all $w_k$ such that $\lambda_k = \lambda_i$.
Therefore,
\begin{align}
      \displaystyle & {\cal G}_i \stackrel{def}{=}
\int_{\R^+}\,\left|\lambda_i(\tilde{w}^2(\tau,z))
-\lambda_i(\tilde{w}^1(\tau,z))\right|\,\md \tau \nonumber \\
\displaystyle \leq &\; \sum_{\tilde{\lambda}_k\neq 0}^n\;\int_{\R^+}\;\left|\frac{\pa \lambda_i} {\pa w_k}\right|\cdot \left|w_k^{2,0}(X^{2,0}(z-\tilde{\lambda}_k\,\tau))-w_k^{1,0}(X^{1,0}(z-\tilde{\lambda}_k\,\tau))\right|\;\md \tau \nonumber \\
 \leq &\; C\;\sum_{\tilde{\lambda}_k\neq 0}^n\;\int_{\R^+}
\; \left|w_k^{2,0}(X^{2,0}(z-\tilde{\lambda}_k\,\tau))
-w_k^{1,0}(X^{1,0}(z-\tilde{\lambda}_k\,\tau))\right|\;\md \tau. \nonumber
\end{align}
Hence, it is still possible to make the change of variables
$\tau \longmapsto z-\tilde{\lambda}_k\,\tau$ to obtain ${\cal G}_i \longrightarrow 0$
as $R_0 \rightarrow 0$. This proves (\ref{5a}) for $X^1$ and $X^2$.

Finally, similarly to the proof of Lemma 4.2 in \cite{PY10}, noting
that (\ref{5.0}) and (\ref{5.2}), we can apply Lemma \ref{L4.2} to
get the desired (\ref{conver'}). Thus, the proof of Lemma \ref{L5.3}
is completed. \hfill $\Box$

\begin{lemma}
\label{L5.4}
Under the assumptions of Proposition \ref{P1.2}, we define
$$  \forall \;(t,x)\in \R^+\times\R,\quad g_t(x)
=X^{1,0}(Z^1(t,X^2(t,Z^{2,0}(x)+\tilde{\lambda}_i\,t))-\tilde{\lambda}_i\,t).$$
Then there is a constant $\gamma>0$ such that
\stepcounter{exte}
\begin{equation}
\label{5g}
\text{a.e.}\;(t,x)\in \R^+\times\R,\quad \frac{\md g_t(x)}{\md x} \geq \gamma
\quad \text{and} \quad g_t(x)\longrightarrow x,\quad \text{as}\quad R_0\rightarrow 0,
\end{equation}
where the convergence is uniform on $\R^+\times\R$.
If furthermore, $w^{j,0}\in W^{1,1}(\R)$ for $j=1,2$, then
\stepcounter{exte}
\begin{equation}
\label{strogt}
    \big|g_t(x) -x\big| \leq C R_0, \quad \forall \, t \geq 0, \; \forall \, x \in \R.
\end{equation}
\end{lemma}

\noindent {\bf Proof.} The proof is omitted here since it is similar to
that of Lemma 4.3 in \cite{PY10}. \hfill $\Box$

\vspace{3mm}

\noindent {\bf Proof of Theorem \ref{T4.1}.}
It suffices to establish the $L^1$ stability for each $w_i\;(i\in K_n)$.
From the explicit expression of the entropy solution, we have
\stepcounter{exte}
\begin{align}
\nonumber
&\|w_i^2(t,\cdot)-w_i^1(t,\cdot)\|_{L^1(\R)} \\\nonumber
=\;&\,\int_{\R}\left|w_i^{2,0}\big(X^{2,0}(Z^2(t,x)-\tilde{\lambda}_i\,t)\big)
-w_i^{1,0}\big(X^{1,0}(Z^1(t,x)-\tilde{\lambda}_i\,t)\big)\right|\,\md
x\\\nonumber \leq \;&
\,\int_{\R}\left|w_i^{2,0}\big(X^{2,0}(Z^2(t,x)-\tilde{\lambda}_i\,t)\big)
-w_i^{1,0}\big(X^{2,0}(Z^2(t,x)-\tilde{\lambda}_i\,t)\big)\right|\,\md
x\\\label{w12i} & \, +
\,\int_{\R}\left|w_i^{1,0}\big(X^{2,0}(Z^2(t,x)-\tilde{\lambda}_i\,t)\big)
-w_i^{1,0}\big(X^{1,0}(Z^1(t,x)-\tilde{\lambda}_i\,t)\big)\right|\,\md
x.
\end{align}
Using the change of variable
$$    \eta_i=X^{2,0}(Z^2(t,x)-\tilde{\lambda}_i\,t),     $$
we get
$$x=X^2(t,Z^{2,0}(\eta_i)+\tilde{\lambda}_i\,t)\quad \text{and}\quad \md x
=\bar{J}^{2,0}(t,\eta_i)\,\md \eta_i,$$
where
\stepcounter{exte}
\begin{equation}
\label{5l}
\bar{J}^{2,0}(t,\eta_i)=\dfrac{N(w^{2,0}(\eta_i))}{N(w^{2,0}(X^{2,0}(Z^{2,0}(\eta_i)
+\tilde{\lambda}_i\,t)))}\geq \gamma > 0,
\end{equation}
with $\gamma$ being a constant. It is easy to see that
\stepcounter{exte}
\begin{align}
\nonumber
&
\,\int_{\R}\left|w_i^{2,0}\big(X^{2,0}(Z^2(t,x)-\tilde{\lambda}_i\,t)\big)
-w_i^{1,0}\big(X^{2,0}(Z^2(t,x)-\tilde{\lambda}_i\,t)\big)\right|\,\md
x\\\nonumber
=\;& \int_{\R} \left|w_i^{2,0}(\eta_i)
-w_i^{1,0}(\eta_i)\right| \bar{J}^{2,0}(t,\eta_i)\,\md \eta_i\\\label{5h}
\leq \;& C\; \|w_i^{2,0}-w_i^{1,0}\|_{L^1(\R)}.
\end{align}
For the second term of the right hand side of (\ref{w12i}), we have
\stepcounter{exte}
\begin{align}
\nonumber
&\;\int_{\R}\left|w_i^{1,0}\big(X^{2,0}(Z^2(t,x)-\tilde{\lambda}_i\,t)\big)
-w_i^{1,0}\big(X^{1,0}(Z^1(t,x)-\tilde{\lambda}_i\,t)\big)\right|\,\md
x\\\label{5.3}
=\;& \int_{\R}\left|w_i^{1,0}(\eta_i)
-w_i^{1,0}(g_t(\eta_i))\right|\bar{J}^{2,0}(t,\eta_i)\,\md \eta_i.
\end{align}
Together with $w_i^{1,0}\in L^1(\R)$ and (\ref{5g}), we may apply
Lemma \ref{L3.7} to get
$$\int_{\R}\left|w_i^{1,0}\big(X^{2,0}(Z^2(t,x)-\tilde{\lambda}_i\,t)\big)
-w_i^{1,0}\big(X^{1,0}(Z^1(t,x)-\tilde{\lambda}_i\,t)\big)\right|\,\md
x\longrightarrow 0,\quad \text{as} \quad R_0\rightarrow 0.$$ This
shows that
$$\|w_i^2(t,\cdot)-w_i^1(t,\cdot)\|_{L^1(\R)}\longrightarrow 0,\quad \text{as}
\quad R_0\rightarrow 0.$$

If $w^{j,0}\in W^{1,1}(\R)$ for $j=1,2$, with the aid of
(\ref{strogt}), we may apply Lemma \ref{L4.2} to (\ref{5.3}). Thus,
together with (\ref{w12i})-(\ref{5.3}), it gives (\ref{stability'}).
\hfill $\Box$

%$5%%%%%%%%%%%%%%%%%%%%%%%%%%%%%%%%%%%%%%%%%%%%%%%%%%%%%%%%%%%%%%%%%%%%%%

\vspace{5mm}

\section{Examples}
\newcounter{exam}
\renewcommand{\theequation}{\thesection.\theexam}

We present two examples for the results obtained in sections 2-4. They
are the Born-Infeld system (of 4 equations) and the augmented Born-Infeld
system (of 8 equations) introduced by Brenier \cite{Br04}, see also
\cite{Se04}. We mention that both the Born-Infeld system and the augmented
Born-Infeld system are rich and linearly degenerate with constant multiplicity
eigenvalues. Moreover, they have a common structure~: $N_i=N$ for all
$i \in K_n$ (see \cite{LPR09}). Another example is the generalized extremal
surface equations, which have been analyzed in detail in \cite{PY10}.

%$5.1%%%%%%%%%%%%%%%%%%%%%%%%%%%%%%%%%%%%%%%%%%%%%%%%%%%%%%%%%%%%%%%%%%%%%%

\subsection{The Born-Infeld system}

The one-dimensional Born-Infeld system reads (see \cite{Br04})~:

\stepcounter{exam}
\begin{equation}
\label{5n}
\left\{\begin{array}{lllll}
 \dis \pa_t D_2 + \pa_x \left(\frac{B_3+D_2P_1-D_1P_2}{h}\right) = 0, \\[4mm]
 \displaystyle{\pa_t D_3 + \pa_x \left(\frac{-B_2+D_3P_1-D_1P_3}{h}\right)
= 0,}\\[4mm]
 \dis \pa_t B_2 + \pa_x \left(\frac{-D_3+B_2P_1-B_1P_2}{h}\right) = 0, \\[4mm]
 \dis \pa_t B_3 + \pa_x \left(\frac{D_2+B_3P_1-B_1P_3}{h}\right) = 0, \\[4mm]
 \dis P(u)=D\times B,\quad h(u)=\sqrt{1+|B|^2+|D|^2+|D\times B|^2}.
\end{array}\right.
\end{equation}
Here $u=(D_2,D_3,B_2,B_3)^{\top}$ are the unknown variables, $B_1$, $D_1$ are real
constants, and
$$B=(B_1,B_2,B_3)^{\top},\quad  D=(D_1,D_2,D_3)^{\top},\quad P=(P_1,P_2,P_3)^{\top}.$$
From \cite{Pe07},
system (\ref{5n}) can be written in the form (\ref{diag}). It is not strictly
hyperbolic with two eigenvalues, denoted by $\lambda$ and $\mu$, of both constant
multiplicity $2$. The variables $\lambda$ and $\mu$ satisfy
\stepcounter{exam}
\begin{equation}
\label{5s}
\left\{\begin{array}{ll}
 \dis    \pa_t \mu + \lambda\; \pa_x \mu = 0, \\[2mm]
 \dis \pa_t \lambda + \mu\; \pa_x \lambda = 0.
 \end{array}\right.
\end{equation}
Moreover, system (\ref{5n}) admits an additional conservation law (see \cite{Br04})~:
$$  \pa_t h + \pa_x P_1=0.   $$
From \cite{LPR09}, (\ref{5n}) is a linearly degenerate rich system with
$$ M=P_1, \quad N_i=N \stackrel{def}{=} h = \frac{2a}{\mu-\lambda}, \quad \mu>\lambda. $$
Here $N_i$ is independent of $i$ for all $i \in K_n$. The eigenvalues
in Lagrangian coordinates are
\stepcounter{exam}
\begin{equation}
\label{F5.3}
 \tilde{\lambda}_1 = \tilde{\lambda}_2 = - a, \quad \tilde{\lambda}_3
= \tilde{\lambda}_4 = a,
\end{equation}
where $a= \sqrt{1+B_1^2+D_1^2}$.

Theorems \ref{T3.1} and \ref{T4.1} can be applied to the Born-Infeld system.
More precisely, given an initial condition for (\ref{5n})~:
\stepcounter{exam}
\begin{equation}
\label{5w}
     t=0: \quad u=u^0(x), \quad x \in \R,
\end{equation}
with $u^0 - \bar{u} \in L^1(\R)\cap L^{\infty}(\R)$ for some constant
$\bar{u} \in \R^4$. We infer from Theorem \ref{T4.1} that the unique entropy
solution of (\ref{5n}) and (\ref{5w}) is $L^1$ stable. Moreover, if
$u^0- \bar{u} \in W^{1,1}(\R)$, we get the strong stability inequality like
(\ref{rtr1}). In addition, when the initial data satisfy the assumptions of
Theorem \ref{T3.1}, the explicit asymptotic limit can be obtained. In this
particular case, we can present better characterization of the explicit
asymptotic limit.

For this purpose, we denote $\lambda^0=\lambda(u^0)$ and $\mu^0=\mu(u^0)$
and suppose
\stepcounter{exam}
\begin{equation}
\label{F5.5}
  \inf_{x \in \R} \mu^0(x) > \sup_{x \in \R} \lambda^0(x).
\end{equation}
Define
$$   Z^0(x) \overset{\text{def}}{=} \int_0^x h(u^0(\xi))\,\md \xi=\int_0^x
\frac{2a}{\mu^0(\xi)-\lambda^0(\xi)}\,\md \xi    $$
and $X^0$ being the inverse function of $Z^0$. It follows from (\ref{X2}) and
$X(0,z) = X^0(z)$ that
\stepcounter{exam}
\begin{equation}
\label{5z}
\dis X(t,z)=\frac{1}{2a}\,\int_0^{z+a\,t}\,\mu^0(X^0(\xi))\,\md \xi
-\frac{1}{2a}\,\int_0^{z-a\,t}\,\lambda^0(X^0(\xi))\,\md \xi.
\end{equation}
This formula was first obtained in \cite{Pe98}. Suppose now
\stepcounter{exam}
\begin{equation}
\label{5aa}
    \dis \lim_{x\rightarrow +\infty}\mu^0(x)=\bar{\mu}^0_+,
\quad \lim_{x\rightarrow -\infty}\lambda^0(x)=\bar{\lambda}^0_-,
\end{equation}
where $\bar{\mu}^0_+ = \mu(\bar{u})$ and $\bar{\lambda}^0_-=\lambda(\bar{u})$.

\begin{remark}
  Assumption (\ref{5aa}) requires the limit only at one direction for each
variable. It is weaker than $\dis{\lim_{x \rightarrow \pm \infty} w^0(x) = \bar{w}}$,
which is made in Theorem \ref{T3.1} and Lemmas \ref{L3.3}-\ref{L3.6}.
\end{remark}

Let $H$ and $I$ be the primitives of functions
$y\longmapsto \lambda^0(X^0(y))-\bar{\lambda}^0_-$ and
$y\longmapsto \mu^0(X^0(y))-\bar{\mu}^0_+$, defined by
\stepcounter{exam}
\begin{equation}
\label{5ab}
 \dis H(y)=\int_1^y\,\big(\lambda^0(X^0(\xi))-\bar{\lambda}^0_-\big)\,\md \xi,
\quad I(y) =\int_1^y\,\big(\mu^0(X^0(\xi))-\bar{\mu}^0_+\big)\,\md \xi,
\quad\forall \;  y\;\in \R.
\end{equation}
Since $H' \in L^1(\R)$ and $I' \in L^1(\R)$, there are two
constants $\bar{H}^-$ and $\bar{I}^+$ such that
\stepcounter{exam}
\begin{equation}
 \label{5ac}
 \dis \lim_{y\rightarrow -\infty}H(y)=\bar{H}^-,
\quad \lim_{y\rightarrow +\infty}I(y)=\bar{I}^+.
\end{equation}
Define
\stepcounter{exam}
\begin{equation}
 \label{5ad}
 \dis \psi(x)=x+\frac{1}{2a}\Big(H(Z^0(x))-\bar{H}^-\Big),\quad \varphi(x)
=x-\frac{1}{2a}\Big(I(Z^0(x))-\bar{I}^+\Big),\quad\forall \;  x\;\in \R.
\end{equation}
A straightforward computation shows that
$$\dis \psi'(x)=\frac{\mu^0(x)-\bar{\lambda}^0_-}{\mu^0(x)-\lambda^0(x)}>0,
\quad \varphi'(x)
=\frac{\bar{\mu}^0_+-\lambda^0(x)}{\mu^0(x)-\lambda^0(x)}>0.   $$
Therefore, $\psi$ and $\varphi$ are all Lipschitz and strictly increasing functions
on $\R$ without any smallness condition on $u^0$.
Let us denote by $\phi$ and $\chi$ be the inverse functions of $\psi$ and $\varphi$,
respectively, which are both Lipschitzian and strictly increasing on $\R$.

It follows from (\ref{5z}) and (\ref{5ac})-(\ref{5ad}) that, for all $x\in \R$,
\begin{align}
\nonumber
&\;X(t,Z^0(x)-a\,t)-\bar{\lambda}^0_-\,t \\\nonumber
=\;& \frac{1}{2a}\int_0^{Z^0(x)}\, \mu^0(X^0(\xi))\,\md \xi-\frac{1}{2a}\int_0^{Z^0(x)-2a\,t}
\, \lambda^0(X^0(\xi))\,\md \xi-\bar{\lambda}^0_-\,t \\\nonumber
=\;& x+\frac{1}{2a}\Big(H(Z^0(x))-H(Z^0(x)-2a\,t)\Big)\\\nonumber
=\;& \psi(x)-\frac{1}{2a}\Big(H(Z^0(x)-2a\,t)-\bar{H}^-\Big)\\\nonumber
\longrightarrow \;& \psi(x),\quad \text{as}\quad t\rightarrow +\infty,
\end{align}
which implies that
\stepcounter{exam}
\begin{equation}
 \label{5ae}
 \dis \phi(X(t,Z^0(x)-a\,t)-\bar{\lambda}^0_-\,t) \longrightarrow x,
\quad \text{as}\quad t\rightarrow +\infty,
\end{equation}
pointwisely on $\R$ and uniformly on any compact set of $\R$.
Similarly,
$$ \dis \chi(X(t,Z^0(x)+a\,t)-\bar{\mu}^0_+\,t) \longrightarrow x,
\quad \text{as}\quad t\rightarrow +\infty,  $$
pointwisely on $\R$ and uniformly on any compact set of $\R$. Let $w_1$ and
$w_2$ be the Riemann invariants associated to the eigenvalue $\lambda$, $w_3$ and
$w_4$ be those associated to $\mu$ (see \cite{Pe07} for the detailed expressions).
Then the asymptotic limit of $w_i \,(1 \leq i \leq 4)$ is
\stepcounter{exam}
\begin{equation}
\label{5af1}
\left\{\begin{array}{ll}
\dis \|w_j(t,\cdot)-w_j^0(\phi(\cdot-\bar{\lambda}^0_-\,t))\|_{L^1(\R)}\longrightarrow 0,
\quad \text{as}\quad t\rightarrow +\infty,  \\[3mm]
\dis \|w_{j+2}(t,\cdot)-w_{j+2}^0(\chi(\cdot-\bar{\mu}^0_+\,t))\|_{L^1(\R)}\longrightarrow 0,
\quad \text{as}\quad t\rightarrow +\infty,
 \end{array}\right.
\end{equation}
for $j=1,2$.

From (\ref{F5.3}), no eigenvalue of the Born-Infeld system in Lagrangian
coordinates is zero. Applying Theorem \ref{T3.1} to (\ref{5n}) together with
the discussion above, we obtain the following result.

\begin{corol}
    Let $u^0-\bar{u} \in L^1(\R) \cap L^\infty(\R)$. Assume (\ref{F5.5}) and (\ref{5aa})
hold. Then the entropy solution of (\ref{5n}) and (\ref{5w}) satisfies (\ref{5af1}).
\end{corol}

%$5.2%%%%%%%%%%%%%%%%%%%%%%%%%%%%%%%%%%%%%%%%%%%%%%%%%%%%%%%%%%%%%%%%%%%%%%

\subsection{The augmented Born-Infeld system}

The augmented Born-Infeld system reads (see \cite{Br04})~:
\stepcounter{exam}
\begin{equation}
\label{5ag}
\left\{\begin{array}{llllllll}
 \dis \pa_t h + \pa_x P_1 = 0, \\[4mm]
 \dis \pa_t P_1 + \pa_x \left(\frac{P_1^2-a^2}{h}\right) = 0, \\[4mm]
 \dis \pa_t D_2 + \pa_x \left(\frac{B_3+D_2P_1-D_1P_2}{h}\right) = 0, \\[4mm]
 \displaystyle{\pa_t D_3 + \pa_x \left(\frac{-B_2+D_3P_1-D_1P_3}{h}\right) = 0,}\\[4mm]
 \dis \pa_t B_2 + \pa_x \left(\frac{-D_3+B_2P_1-B_1P_2}{h}\right) = 0, \\[4mm]
 \dis \pa_t B_3 + \pa_x \left(\frac{D_2+B_3P_1-B_1P_3}{h}\right) = 0, \\[4mm]
 \dis \pa_t P_2 + \pa_x \left(\frac{P_1P_2-D_1D_2-B_1B_2}{h}\right) = 0, \\[4mm]
 \dis \pa_t P_3 + \pa_x \left(\frac{P_1P_3-D_1D_3-B_1B_3}{h}\right) = 0,
\end{array}\right.
\end{equation}
where $h,P_1,D_2,D_3,B_2,B_3,P_2$ and $P_3$ are the unknown variables.
From \cite{Pe07}, the eigenvalues (\ref{5ag}) are~:
\begin{equation}
 \nonumber
 \dis \lambda\overset{\text{def}}{=}\lambda_1=\lambda_2=\lambda_3
=\frac{P_1-a}{h},\; \nu\overset{\text{def}}{=}\lambda_4=\lambda_5
=\frac{P_1}{h}, \; \mu\overset{\text{def}}{=}\lambda_6=\lambda_7
=\lambda_8=\frac{P_1+a}{h}.
\end{equation}
It is easy to see that $\lambda$ and $\mu$ are also the Riemann
invariants satisfying (\ref{5s}). Moreover, the eigenvalues of
system (\ref{5ag}) in Lagrangian coordinates are~:
$$\tilde{\lambda}_1 = \tilde{\lambda}_2 = \tilde{\lambda}_3
= - a, \quad \tilde{\lambda}_4 = \tilde{\lambda}_5
= 0, \quad \tilde{\lambda}_6 = \tilde{\lambda}_7 = \tilde{\lambda}_8 = a.$$

Since the augmented Born-Infeld system and the Born-Infeld system
have a common structure~: $N_i= h\,(1 \leq i \leq 8)$ and the same
eigenvalues $\lambda$ and $\mu$, it suffices to consider the case
for $\tilde{\lambda}_4 = \tilde{\lambda}_5 = 0$. According to the discussion
above, we only need to check that the expression like
(\ref{5ae}) holds for $\nu$.

Indeed, with the notations of the previous subsection, we have
\begin{align}
\nonumber
&\;X(t,Z^0(x))-\frac{\bar{\lambda}^0_-+\bar{\mu}^0_+}{2}\,t \\\nonumber
=\;& \frac{1}{2a}\int_0^{Z^0(x)+a\,t}\, \mu^0(X^0(\xi))\,\md \xi
-\frac{1}{2a}\int_0^{Z^0(x)-a\,t}\, \lambda^0(X^0(\xi))\,\md \xi
-\frac{\bar{\lambda}^0_-+\bar{\mu}^0_+}{2}\,t \\\nonumber
=\;& x+\frac{1}{2a}\Big(I(Z^0(x)+a\,t)-I(Z^0(x))
+\frac{1}{2a}\Big(H(Z^0(x))-H(Z^0(x)-a\,t)\Big)\\\nonumber
\longrightarrow \;& x+\frac{1}{2a}\Big(\big(H(Z^0(x))-\bar{H}^-\big)
-\big(I(Z^0(x))-\bar{I}^+\big)\Big)\overset{\text{def}}{=}\omega(x),\quad
\text{as}\quad t\rightarrow +\infty,
\end{align}
where we have used the
assumption (\ref{5ac}). A straightforward calculation shows that
$$\dis \omega'(x)=\frac{\bar{\mu}_+^0-\bar{\lambda}_-^0}{\mu^0(x)-\lambda^0(x)}>0.$$
Then we can define its inverse function denoted by $\theta$, which is
also a Lipschitz and strictly increasing function on $\R$. Moreover,
$$\dis \theta\left(X(t,Z^0(x))-\frac{\bar{\lambda}^0_-
+\bar{\mu}^0_+}{2}\,t\right)\longrightarrow  x,
\quad \text{as}\quad t\rightarrow +\infty,$$
pointwisely on $\R$ and uniformly on any compact set of $\R$.

Thus, for the eigenvalue $\nu$ which does not appear in the Born-Infeld system,
there still exists a Lipschitz function $\theta$ such that the expression like
(\ref{5ae}) holds.
This shows that for initial data in $L^1(\R)\cap L^{\infty}(\R)$ the entropy solution
of the augmented Born-Infeld system has analogous asymptotic limits to (\ref{5af1}),
without any smallness condition on the initial data. The $L^1$ stability of solutions
can be obtained in a same way as for the Born-Infeld system.

% Appendix ---------------------------------------------------------

%\section*{Acknowledgements}
\vspace{5mm}

\noindent {\large \bf Acknowledgements}

\vspace{2mm}

This work was in part carried out when Yong-Fu YANG had a postdoctor position of the academic year 2009-2010 at Universit\'{e} Blaise Pascal. He would like to thank Laboratoire de Math\'ematiques (CNRS-UMR 6620) for hospitality. His research was also supported in part by National Natural
Science Foundation of China under Grant No. 10926162, Fundamental
Research Funds for the Central Universities No. 2009B01314
and Natural Science
Foundation of Hohai University under Grant No. 2009428011.

% ---------------------------------------------------------

\end{document}